\documentclass[twoside, 11pt]{article}

\usepackage{amssymb, amsmath, latexsym, mathrsfs, verbatim, calc}

\numberwithin{equation}{section} 
\newtheorem{thm}{Theorem}[section]
\newtheorem{lem}[thm]{Lemma}

\newtheorem{prop}[thm]{Proposition}
\newtheorem{rem}[thm]{Remark}

\newtheorem{eg}[thm]{Example}
\newtheorem{defn}[thm]{Definition}

\def\ba{\begin{array}}
\def\ea{\end{array}}
\def\beq{\begin{equation}}
\def\bes{\begin{equation*}}
\def\ees{\end{equation*}}
\def\bea{\begin{eqnarray}}
\def\eea{\end{eqnarray}}
\def\beas{\begin{eqnarray*}}
\def\eeas{\end{eqnarray*}}

\def\dis{\displaystyle}

\def\no{\noindent}

\def\norm{|\nts|}

\def\lastline{\par \vspace{-7.3ex} \no}
\def\restart{\par \vspace{-3.65ex} \no}
\def\nts{\negthinspace}
\def\ss{\smallskip}
\def\ms{\medskip}
\def\bs{\bigskip}
\def\q{\quad}
\def\qq{\qquad}

\def\ua{\mathop{\uparrow}}

\def\={=\nts \nts=\nts \nts=\nts \nts=}

\def\wt{\widetilde}

\def\({\textnormal{(}}
\def\){\textnormal{)}}

\def\cd{\cdot}
\def\cds{\cdots}

\def\fa{\,\forall \,}
\def\pa{\partial}
\def\es{\emptyset}

\def\a{\alpha}
\def\g{\gamma}
\def\d{\delta}

\def\z{\zeta}
\def\k{\kappa}
\def\l{\lambda}
\def\m{\mu}

\def\si{\sigma}
\def\t{\tau}
\def\f{\varphi}
\def\th{\theta}
\def\o{\omega}

\def\f{\phi}

\def\D{\Delta}

\def\L{\Lambda}
\def\O{\Omega}
\def\F{\Phi}
\def\P{\Psi}

\def\bF{{\bf F}}

\def\cD{{\cal D}}
\def\cE{{\cal E}}
\def\cF{{\cal F}}
\def\cG{{\cal G}}
\def\cH{{\cal H}}

\def\cM{{\cal M}}

\def\cT{{\cal T}}

\def\hB{\mathbb{B}}
\def\hC{\mathbb{C}}
\def\hD{\mathbb{D}}
\def\hE{\mathbb{E}}

\def\hN{\mathbb{N}}

\def\hQ{\mathbb{Q}}
\def\hR{\mathbb{R}}

\def\sB{\mathscr{B}}

\def\sD{\mathscr{D}}
\def\sE{\mathscr{E}}

\def\sL{\mathscr{L}}

\def\sP{\mathscr{P}}

\def\esssup{\mathop{\rm esssup}}

\def\dtp{{\hbox{$dt \times dP$-a.s.}}}
\def\pas{{\hbox{$P$-a.s.}}}

\def\no{\noindent}

\def\ss{\smallskip}
\def\ms{\medskip}
\def\bs{\bigskip}
\def\q{\quad}
\def\qq{\qquad}
\def\hb{\hbox}
\def\pa{\partial}
\def\cd{\cdot}
\def\cds{\cdots}

\def\lan{\langle}
\def\ran{\rangle}

\def\bF{{\bf F}}

\def\neg{\negthinspace}
\def\dneg{\neg \neg}
\def\tneg{\dneg \neg}
\def\qed{\hfill \rule[0cm]{.25cm}{.25cm}\medskip}   
\def\dfnn{\stackrel{\triangle}{=}}

\def\lan{\langle}
\def\ran{\rangle}
\def\b1{{\bf 1}}
\def\cad{c\`{a}dl\`{a}g~}

\newenvironment{proof}
{\no {\bf Proof.\;}}{$\qed$ \vspace{-0.2 in}\\}

\newenvironment{itm}{\vspace{-1ex}\begin{itemize}}{\end{itemize}}
\def\bi{\begin{itm}}
\def\ei{\end{itm}}

\def\equ_ind{\arabic{section}.\arabic{equation}}
\def\sec_ind{\arabic{section}}

\begin{document}

\title{\bf Representation Theorems for Quadratic $\cF$-Consistent Nonlinear Expectations
\footnote{\no Part of this work was completed while the second and
third authors were visiting IRMAR,
Universit\'e Rennes 1, France, whose hospitality is greatly
appreciated. }}

\author{
Ying Hu,\thanks{ \no IRMAR, Universit\'e Rennes 1, Campus de Beaulieu,
35042 Rennes Cedex, France;  email: ying.hu@univ-rennes1.fr.}\q
Jin Ma,\thanks{ \noindent Department of Mathematics, Purdue
University, West Lafayette, IN 47907-1395; email:
majin@math.purdue.edu. This author is supported in part by NSF
grant \#0505427. }\q Shige Peng,\thanks{ \no Institute of
Mathematics, Shandong University, 250100, Jinan, China; email:
peng@sdu.edu.cn.} ~~ and ~~ Song Yao\thanks{ \noindent Department
of Mathematics, Purdue University, West Lafayette, IN 47907-1395;
email: syao@purdue.edu. } }

\date{April 12, 2007}

\maketitle

\centerline{\bf Abstract }

\bs

In this paper we extend the notion of ``filtration-consistent
nonlinear expectation" (or ``$\cF$-consistent nonlinear
expectation") to the case when it is allowed to be dominated by a
$g$-expectation that may have a quadratic growth. We show that for
such a nonlinear expectation many fundamental properties of a
martingale can still make sense, including the Doob-Meyer type
decomposition theorem and the optional sampling theorem. More
importantly, we show that
any quadratic $\cF$-consistent nonlinear expectation with a certain
domination property must be a quadratic $g$-expectation as was
studied in Ma-Yao \cite{MY}. The main contribution of this paper is
the finding of the domination condition to replace the one used in
all the previous works (e.g., \cite{CHMP} and \cite{Pln}), which is
no longer valid in the quadratic case. We also show that the
representation generator must be deterministic, continuous, and
actually must be of the simple form $g(z) = \m(1+|z|)|z|$, for some
constant $\m>0$.

\vfill \vspace{.7cm} \no {\bf Keywords: }\: Backward SDEs,
$g$-expectation, $\cF$-consistent nonlinear expectations,
quadratic nonlinear expectations, BMO, Doob-Meyer decomposition,
representation theorem.

\eject

\section{Introduction}
\setcounter{equation}{0}

In this paper we study a class of {\it filtration-consistent
nonlinear expectations} (or {\it $\cF$-consistent nonlinear
expectations}), first introduced in Coquet-Hu-M\'emin-Peng
\cite{CHMP}. Such nonlinear expectations are natural extensions of
the so-called $g$-expectation, initiated in \cite{Peng-97}, and
therefore have direct relations with a fairly large class of risk
measures in finance.  The main point of interest of this paper is
that the nonlinear expectations are allowed to have possible
quadratic growth, and our ultimate goal is to prove a
representation theorem that characterizes the nonlinear
expectations in terms of a class of quadratic BSDEs. We should
note that the class of ``quadratic nonlinear expectations" under
consideration contains many convex risk measures that are not
necessarily ``coherent". The most notable example is the entropic
risk measure (see, e.g., Barrieu and El Karoui \cite{BarKar}),
which is known to have a representation as the solution to a
quadratic BSDE, but falls outside the existing theory of
$\cF$-consistent nonlinear expectations. We refer to \cite{ADEH}
and \cite{FS-02} for the basic concepts of coherent and convex
risk measures, respectively, to
 \cite{Pln}  for detailed accounts of the
relationship between the risk measures and non-linear
expectations. A brief review of the basic properties of
$\cF$-consistent nonlinear expectations will be given in \S2 for
ready references.

An interesting result so far in the development of the notion of
$\cF$-consistent nonlinear expectations is its relationship with the
backward stochastic differential equations (BSDEs). Although as an
extension of the so-called {\it $g$-expectation},
which is defined directly via the BSDE, it is conceivable that an
$\cF$-consistent non-linear expectation should have some
connection to BSDEs, its proof is by no means easy. In the case
when $g$ has only linear growth, it was shown in \cite{CHMP} that
if an $\cF$-consistent non-linear expectation is ``{\it
dominated}" by a $g^\m$-expectation in the sense that
 \bea
 \label{domin}
 \cE[X]-\cE[Y] \le \cE^{g^\m}[X-Y],\qq  \fa X,Y \in L^2(\cF_T)
 \eea
where $g^\m=\mu|z|$ for some constant $\mu>0$,  then it has to be
a $g$-expectation.
The significance of such a result might be more clearly seen from
the following consequence in finance: {\it any coherent risk
measure satisfying the required domination condition can be
represented by the solution of a simple BSDE}(!). In an
accompanying paper by Ma and Yao \cite{MY}, the notion of
$g$-expectation was generalized to the quadratic case, along with
some elementary properties of the $g$-expectations including the
Doob-Meyer decomposition and upcrossing inequalities. However, the
representation property for general (even convex) risk measure
seems to be much more subtle. One of the immediate obstacles is
that the ``domination" condition (\ref{domin}) breaks down in the
quadratic case. For example, one can check that a quadratic $g$
expectation with $g=\m(|z|+|z|^2)$ {\it cannot} be dominated by
itself(!). Therefore some new ideas for replacing the domination
condition (\ref{domin}) are in order.

The main purpose of this paper is to generalize the notion of
$\cF$-consistent nonlinear expectation to quadratic case and prove
at least a version of the representation result for such nonlinear
expectations. An important contribution of this paper is the finding
of a new domination condition for the quadratic nonlinear
expectation, stemmed from the Reverse H\"older Inequality in BMO
theory \cite{Ka}. More precisely, we observe that there exists an
$L^p$ estimation for the difference of quadratic $g$-expectations by
using the reverse H\"older inequality. Extending such an estimate to
the general nonlinear expectations, we then obtain
an $L^p$-type domination which turns out to be sufficient for our
purpose.
Following the idea in \cite{Pln}, with the help of the new
domination condition, we then prove the optional sampling, and a
Doob-Meyer type decomposition theorem for quadratic
$\cF$-martingales.
Similar to the linear case, we can then prove that the
representation property for the quadratic $\cF$-consistent nonlinear
expectation remains valid under such domination condition. That is,
one can always find a quadratic $g$-expectation with $g$ being of
the form: $g=\m(1+|z|)|z|$, to represent the given nonlinear
expectation.

Our discussion on quadratic nonlinear expectation benefited greatly
from the recent development on the theory of BSDEs with quadratic
growth, initiated by Kobylanski \cite{Ko} and the subsequent results
on such BSDEs with unbounded terminal conditions by Briand and Hu
\cite{BH-06,BH-07}. In particular, we need to identify an
appropriate subset of exponentially integrable random variables with
certain algebraic properties on which a quadratic $\cF$-consistent
nonlinear expectation can be defined.
It is worth noting that such a set will have to contain all the
random variables of the form $\xi+z B_\t$, where $B$ is the driving
Brownian motion, $\xi\in L^\infty(\cF_T)$, and $\t$ is any stopping
time, which turns out to be crucial in proving the representation
theorem and the continuity of the representation function $g$.
We should remark that although most of the steps towards our final
result look quite similar to the linear growth case, some special
treatments are necessary along the way to overcome various
technical subtleties caused by the quadratic BSDEs, especially
those with unbounded terminal conditions. We believe that many of
the results are interesting in their own right. We therefore
present full details for future references.



\ms

This paper is organized as follows. In section 2 we give the
preliminaries and review some basics of quadratic $g$-expectations
and the BMO martingales. In section 3 we introduce the notion of
quadratic $\cF$-consistent nonlinear expectations and several new
notions of the dominations. In section 4, we show some properties
of quadratic $\cF$-expectations including the optional sampling
theorem, which pave the ways for the later discussions. In section
5, we prove a Doob-Meyer type decomposition theorem for quadratic
$\cF$-submartingales. The last section is devoted to the proof of
the representation theorem of the quadratic nonlinear
expectations.


\section{Preliminaries}
\setcounter{equation}{0}

Throughout this paper we consider a filtered, complete probability
space $(\O,\cF, P, \bF)$ on which is defined a $d$-dimensional
Brownian motion $B$. We assume that the filtration $\bF\dfnn
\{\cF_t\}_{t\ge0}$ is generated by the Brownian motion $B$,
augmented by all the $P$-null sets in $\cF$, so that it satisfies
the {\it usual hypotheses} (cf. \cite{Pr-90}).  We denote $\sP$ to
be the progressive measurable $\si$-field on $\O\times [0,T]$; and
$\cM_{0,T}^B$ to be the set of all $\bF$-stopping times $\tau$
such that $0\leq\tau\leq T$, $P$-a.s., where $T>0$ is some fixed
time horizon.

In what follows we fix a finite time horizon $T>0$, and denote $\hE$
to be a generic Euclidean space, whose inner products and norms will
be denoted as the same $\langle\cdot,\cdot\rangle$ and $|\cdot|$,
respectively; and denote $\hB$ to be a generic Banach space with
norm $\|\cd\|$. Moreover, we shall denote $\cG\subseteq\cF$ to be
any sub-$\si$-field, and for any $x\in\hR^d$ and any $r>0$ we denote
$\overline{B}_r(x)$ to be the
closed ball with center $x$ and radius $r$.
 Furthermore, the following spaces of functions will be frequently
used in the sequel. We denote

\begin{itemize}

\item for $0\le p\le\infty$, $L^p(\cG;\hE)$ to be the space of all $\hE$-valued,
$\cG$-measurable random variables $\xi$, with $E(|\xi|^p)<\infty$.
In particular, if $p=0$,
then $L^0(\cG,\hE)$ denotes the space of all $\hE$-valued, $\cG$-measurable random
variables; and if $p=\infty$, then $L^\infty(\cG;\hE)$ denotes the space of all
$\hE$-valued, $\cG$-measurable random variables $\xi$ such that
$\|\xi\|_\infty \dfnn \underset{\o \in \O}{\esssup}|\xi(\o)|<\infty$;

\item $0 \le  p\le\infty$, $L^p_\bF([0,T];\hB)$ to be the space of all
$\hB$-valued, $\bF$-adapted processes $\psi$, such that
$E\int_0^T\|\psi_t\|^pdt<\infty$. In particular, $p=0$ stands for
all $\hB$-valued, $\bF$-adapted processes; and $p=\infty$ denotes
all processes $X\in L^0_\bF([0,T];\hB)$
such that $\|X\|_\infty \dfnn \underset{t,\o}
{\esssup} |X(t,\o)|<\infty$;

\item $\hD^\infty_\bF([0,T];\hB)=\{X\in L^\infty_\bF([0,T];\hB): \hb{ $X$ has c\`adl\`ag paths}\}$;

\item $\hC^\infty_\bF([0,T];\hB)=\{X\in \hD^\infty_\bF([0,T];\hB): \hb{ $X$ has continuous
paths}\}$;

\item $\cH^2_\bF([0,T];\hB)=\{X \in L^2_\bF([0,T];\hB): \hb{$X$ is predictably
measurable}\}$.
\end{itemize}

The following two spaces are variations of the $L^p$ spaces defined above, they will be
important for our discussions regarding quadratic BSDEs with unbounded terminal conditions.
For any $p>0$, we denote $\cM^p(\hR^d)$ to be the space of all $\hR^d$-valued predictable
processes $X$ such that
 \bea
 \label{Mp}
 \|X \|_{\cM^p} \dfnn \Big(E\big(\int_0^T |X_s|^2ds \big)^{p/2}\Big)^{1 \land 1/p}
  < \infty.
 \eea
We note that for $p \geq 1$, $\cM^p(\hR^d)$ 
is a Banach space with the norm  $\|\cd\|_{\cM^p}$, 
and for $p \in (0,1)$, $\cM^p(\hR^d)$ 
is a complete metric space with the distance defined through
(\ref{Mp}) 
Finally,
if $d=1$, we shall drop $\hE=\hR$ from the notation (e.g.,
$L^p_\bF([0,T])=L^p_\bF([0,T];\hR)$,
$L^\infty({\cF_T})=L^\infty({\cF_T};\hR)$, and so on).

\bs

\no{\bf Quadratic $g$-expectations on $L^\infty({\cF_T})$}~

\ms

We now give a brief review of the notion of {\it quadratic
$g$-expectations} studied in Ma and Yao \cite{MY}. First recall that
for any $\xi\in L^2({\cF_T})$, and a given ``generator"
$g=g(t,\o,y,z): [0,T] \times \O \times \hR \times \hR^d \mapsto \hR$
satisfying the standard conditions (e.g., it is Lipschitz in all
spatial variables, and is of linear growth, etc.), the
$g$-expectation of $\xi$ is defined as $\cE^g(\xi)\dfnn Y_0$, where
$Y=\{Y_t:0 \le t \le T \}$ is the solution to the following BSDE:
 \bea
 \label{BSDE}
 Y_t=\xi+\int_t^T g(s,Y_s,Z_s)ds - \int_t^T Z_s dB_s,  \qq \fa t \in
 [0,T].
 \eea
We shall denote (\ref{BSDE}) by BSDE$(\xi,g)$ in the sequel for
notational convenience.

In \cite{MY} the $g$-expectation was extended to the quadratic
case, based on the well-posedness result of the quadratic BSDEs by
Kobylanski \cite{Ko}, and under rather general conditions on the
generator $g$. In this paper, however, we shall be content
ourselves with a slightly simplified form of the generator $g$
that is sufficient for our purpose. More precisely, we assume that
the generator $g$ is independent of the variable $y$, and
satisfies the following {\it Standing Assumptions}:
 \bi
 \item[{\bf (H1)}] The function $g: [0,T]\times\O\times\hR^d\mapsto \hR$ is
$\sP \otimes \sB(\hR^d)$-measurable
 and $g(t,\o,\cdot)$ is continuous for all $(t,\o) \in [0,T] \times \O$;

 \item[{\bf (H2)}] There exists a constant $\ell>0$ such that for $dt \times dP$-a.s.
 $(t,\o) \in [0,T] \times \O$ and any $z \in \hR^d$
 \bea
 \label{H2}
 |g(t,\o,z)|\leq \ell(|z|+|z|^2) \mbox{\q and \q} \Big|\frac{\pa g}{\pa
 z}(t,\o,z)\Big| \leq \ell(1+|z|).
 \eea
 \ei

%
%

In light of the results of \cite{Ko} we know that under the
assumptions (H1) and (H2), for any $\xi\in L^\infty({\cF_T})$ the
BSDE (\ref{BSDE}) has a unique solution $(Y,Z) \in
\hC^\infty_\bF([0,T])\times \cH^2_\bF([0,T];\hR^d)$. We can then
define the quadratic $g$-expectation of $\xi$ as $\cE^g(\xi)=Y_0$
and the {\it conditional
$g$-expectation} as
 \bea
 \label{condigexp}
 \cE^g[\xi|\cF_t] \dfnn Y^\xi_t, \q~~ \fa t \in [0,T], \q \fa \xi
 \in L^\infty(\cF_T).
 \eea

It is easy to see that $g|_{z=0}=0$ from (H2). So by the uniqueness
of the solution to the quadratic BSDE, one can show that all the
fundamental properties of nonlinear expectations are still valid for
quadratic $g$-expectations:
 \bi
 \item[{(i)}]
 (Time-consistency)~~$\cE^g\big[\cE^g[\xi|\cF_t]\big|\cF_s\big]=\cE^g[\xi|\cF_s],
 \q \pas \q \fa 0\le s \le t$;

 \item[{(ii)}] (Constant-preserving)~~$\cE^g[\xi|\cF_t]=\xi,\q \pas
 \q \fa \xi \in L^\infty({\cF_t})$;

 \item[{(iii)}] (`` Zero-one Law")~~$\cE^g[\b1_A \xi |\cF_t]=\b1_A
 \cE^g[ \xi | \cF_t],\q \pas \q \fa A \in \cF_t$.
 \ei
Furthermore, since $g$ is independent of $y$, then we know that
the quadratic $g$-expectation is also ``{\it translation
invariant}" in the sense that
 \bea
 \label{transinv}
 \cE^g[\xi+\eta|\cF_t]=\cE^g[\xi|\cF_t]+\eta,\q \pas ~~~ \fa t \in
 [0,T], ~\fa \eta \in L^\infty({\cF_t}).
 \eea

Along the same lines of \cite{Pln} we can define the ``quadratic
$g$-martingales" as usual. For example,
A process $X \in L^\infty_\bF([0,T])$ is called a
$g$-submartingale (resp. $g$-supermartingale) if for any $0 \leq s
< t \leq T$, it holds that
 \beas
 \cE^g[X_t|\cF_s] \geq ~ (\mbox{resp. }\le)~ X_s, \qq \pas
 \eeas
The process $X$ is called a quadratic $g$-martingale if it is both
a $g$-submartingale and a $g$-supermartingale.

Similar to the cases studied in \cite{Pln} where $g$ is Lipschitz
continuous and of linear growth, it was shown in \cite{MY} that the
quadratic $g$-sub(super)martingales also admit the Doob-Meyer type
decomposition, and an upcrossing inequality holds (cf. \cite[Theorem
4.6]{MY}). The next theorem summarizes some results of \cite{MY},
adapted to the current setting, which will be used in our future
discussion. The proof of these results can be found in \cite[Theorem
4.2 and Corollary 4.7]{MY}.
\begin{thm}
 \label{mayaothm}
 Assume (H1) and (H2). Then, for
 any right-continuous $g$-submartingale  (resp.
$g$-supermartingale) $Y \in L^\infty_\bF([0,T])$, there exist an
increasing (resp. decreasing) \cad process $A$ null at $0$ and a
process $Z \in \cH^2_\bF([0,T];\hR^d)$, such that
 \beas
 Y_t=Y_T+\int_t^T g(s,Y_s,Z_s)ds-A_T+A_t-\int_t^T Z_s dB_s, \qq
 t\in[0,T].
 \eeas
Furthermore, if $g$
vanishes as $z$ vanishes, then any
$g$-submartingale (resp. $g$-supermartingale) $X$ must satisfy the
following continuity property: For any dense subset $\cD$ of
$[0,T]$, $P$-almost surely, the limit $\dis
\underset{r \nearrow t,\, r \in \cD}{\lim} X_r$ (resp. $\underset{r
\searrow t,\, r \in \cD}{\lim} X_r$) exists for any $t \in (0,T]$
(resp. $t \in [0,T)$).
\qed
\end{thm}

\bs

\no{\bf BMO and Exponential Martingales}
\ms

To end this section, we recall some important facts regarding the
so-called ``BMO martingales" and the properties of the related
stochastic exponentials. We refer to the monograph of Kazamaki
\cite{Ka} for a complete exposition of the theory of continuous BMO
and exponential martingales. Here we shall be content with only some
facts that are useful in our future discussions.

To begin with, we recall that a uniformly integrable martingale
$M$ null at zero is called a ``BMO martingale" on $[0,T]$, if for
some $1\le p<\infty$, it holds that
 \bea
 \label{BMOp}
 \|M\|_{BMO_p}\dfnn\sup_{\t\in\cM_{0,T}}\Big\|E\{|M_T-M_{\t-}|^p\big|
 \cF_\t\}^{1/p}\Big\|_\infty<\infty.
 \eea
In such a  case we denote $M\in$BMO$(p)$. It is important to note
that $M\in$BMO$(p)$ if and only if $M\in$BMO$(1)$, and all the
BMO$(p)$ norms are equivalent. Therefore in what follows we shall
say that a martingale $M$ is BMO without specifying the index $p$;
and we shall use only the BMO$(2)$ norm and denote it simply by
$\|\cd\|_{BMO}$. Note also that for a {\it continuous} martingale
$M$ one has
 \beas
 \|M\|_{BMO}=\|M\|_{BMO_2}=\sup_{\t\in\cM_{0,T}}
\Big\|E\{\lan M\ran_T-\lan M\ran_\t\big|\cF_\t\}\Big\|_\infty.
 \eeas

Now, for a given Brownian motion $B$, we say that a process $Z\in
L^2_\bF([0,T];\hR^d)$ is a BMO process, denoted by $Z \in$ BMO
with a slight abuse of notations,
if the stochastic integral $M\dfnn Z\bullet B$ is a BMO martingale.
We remark that the space of BMO martingales is smaller than any
$\cM^p(\hR^d)$ space (see (\ref{Mp}) for definition). To wit, it
holds that BMO$ \subset \underset{p>0}{\bigcap} \cM^p(\hR^d)$.
Furthermore, by the so-called ``Energy Inequality" \cite[p.29]{Ka},
one checks that
 \bea
 \label{BMOdom}
 \big(\|Z\|_{\cM_{2n}}\big)^{2n}=E\Big(\int_0^T |Z_s|^2 ds\Big)^n \leq
 n!\|Z \|^{2n}_{BMO},\q~~ \fa n \in \hN.
 \eea

We now turn our attention to the stochastic exponentials of the
BMO martingales. Recall that for a continuous martingale $M$, the
Dol\'eans-Dade stochastic exponential of $M$, denoted customarily
by $\sE(M)$, is defined as $\sE(M)_t\dfnn \exp\{M_t-\frac12\lan
M\ran_t\}$, $t\ge 0$. Note that if $\sE(M)$ is a uniformly
integrable martingale, then the H\"older inequality implies that
 \bea
 \label{holder}
 \sE(M)^p_\t\leq E[\sE(M)^p_T|\cF_\t],\qq \pas
 \eea
for any stopping time $\t \in \cM_{0,T}$ and any $p \geq 1$.
However, if $M$ is further a BMO martingale, then the stochastic
exponential $\sE(M)$ is itself a uniform integrable martingale (see
\cite[Theorem 2.3]{Ka}). Moreover, the so-called ``Reverse H\"older
Inequality" (cf. \cite[Theorem 3.1]{Ka}) holds for $\sE(M)$. We note
that this inequality plays a fundamental role in the new domination
condition for the nonlinear expectations, which leads to the
representation theorem and its continuity, we give the complete
statement here for ready references.
For any $\a>2$, define
 \bea
 \label{phialpha}
 \f_\a(x) \dfnn \bigg\{1+x^{-2}\log{\Big[(1-2
 \a^{-x})\frac{2x-1}{2x-2}\Big]}\bigg\}^{\frac{1}{2}}-1, \qq x \in
 (1,\infty).
 \eea

\begin{thm}  \mbox{\rm (Reverse H\"older Inequality)}~
Suppose that  $M\in$ BMO$(p)$ for $p \in (1,\infty)$. If it
satisfies that $\| M \|_{BMO} \leq \f_\a(p)$, then one has
 \bea
 \label{rhi}
 E\big[\sE(M)^p_T \big|\cF_\t\big]\leq \a^p \sE(M)^p_\t, \qq \fa
 \t \in \cM^B_{0,T}.
 \eea
\end{thm}

Finally, we give a result that relates the solution to a quadratic
BSDE to the BMO processes. Let us consider the BSDE (\ref{BSDE})
in which the generator $g$ has a quadratic growth. For simplicity,
we assume
there is some $k>0$ (we may assume without loss of generality that
$k \ge \frac{1}{2}$) such that for $dt \times dP$-a.s. $(t,\o) \in
[0,T] \times \O$,
 \bea
 \label{gquad}
  |g(t,\o,y,z)| \leq k(1+|z|^2), \qq \fa (y,z) \in \hR \times \hR^d.
 \eea
%
Let $(Y,Z)\in \hC^\infty_\bF([0,T]) \times \cH^2_\bF([0,T];\hR^d)$
be a solution to (\ref{BSDE}). Applying It\^o's formula to $\dis
e^{4kY_t}$ from $t$ to $T$ one has:
 \beas
 e^{4kY_t}+8k^2\int_t^T e^{4kY_s}|Z_s|^2ds &=& e^{4kY_T} +
 4k\int_t^T e^{4kY_s}g(s,Y_s,Z_s)ds-4k\int_t^T e^{4kY_s}Z_sdB_s
 \nonumber \\
 &\leq& e^{4kY_T}+4k^2\int_t^Te^{4kY_s}\big(1+|Z_s|^2\big)ds
  -4k\int_t^T e^{4kY_s}Z_sdB_s.
  \eeas
Taking the conditional expectation $E\{\cdot|\cF_t\}$ on both
sides above, and then use some standard manipulations one derives
fairly easily that
 \beas
 E\big[\int_t^T |Z_s|^2 ds|\cF_t\big] \leq e^{4k \|Y \|_\infty}
 E\big[e^{4k\xi}-e^{4kY_t}|\cF_t\big]+e^{8k \| Y \|_\infty}(T-t).
 \eeas
In other words, we have proved the following result.
\begin{prop}
\label{BMO1}
 Suppose that $(Y,Z)\in \hC^\infty_\bF([0,T]) \times \cH^2_\bF([0,T];\hR^d)$ is a
 solution of the BSDE (\ref{BSDE})
 with $\xi \in L^\infty({\cF_T})$, and
 $g$ satisfies (\ref{gquad}). Then $Z\in$ BMO, and the following
 estimate holds:
 \beas
 \|Z \|^2_{BMO} \leq (1+T)e^{8k \| Y \|_\infty}.
\eeas
\mbox{} \lastline
\qed
\end{prop}

\section{Quadratic $\cF$-Expectations}
\setcounter{equation}{0}

In this section we introduce the notion of ``quadratic
$\cF$-consistent nonlinear expectation". To begin with, we recall
from \cite{Pln} that an $\cF$-consistent nonlinear expectation is a
family of operators, denoted by $\{\cE_t\}_{t\ge 0}$, such that for
each $t\in[0,T]$, $\cE_t: L^0(\cF_T)\mapsto L^0(\cF_t)$, and that
the following axioms are fulfilled:
 \bi
 \item[{\bf (A1)}] {\bf Monotonicity:}~~$\cE_t[\xi] \geq \cE_t[\eta ]$, $\pas$, if
 $\xi \geq \eta$, $\pas$;
 \item[{\bf (A2)}] {\bf Constant-Preserving:}~~$\cE_t[\xi]=\xi$, $\pas$, $\fa \xi
 \in L^0(\cF_t)$;
 \item[{\bf (A3)}] {\bf Time-Consistency:}~~$\cE_s[\cE_t[\xi]]=\cE_s[\xi]$, $\pas$,
 $\fa s \in [0,t]$;
 \item[{\bf (A4)}] {\bf ``Zero-One Law":}~~$\cE_t[\b1_A\xi]=\b1_A\cE_t[\xi]$, $\pas$,
 $\fa A \in \cF_t$.
 \ei
The operator $\cE_t[\cd]$ has been called the ``nonlinear
conditional expectation", and denoted by $\cE\{\cd|\cF_t\}$ for
obvious reasons. It was worth noting that in all the previous
cases the natural ``domain" of the nonlinear expectation is the
space $L^2(\cF_T)$, thus a nonlinear expectation can be related to
the solution to the BSDEs using the ``classical" theory.

In the quadratic case, however, the situation is quite different. In
particular, if the main concern is the representation theorem where
the quadratic BSDE is inevitable, then the domain of the nonlinear
expectation will become a fundamental issue. For example, due to the
limitation of the well-posedness of a quadratic BSDE, a quadratic
nonlinear expectation would naturally be restricted to the space
$L^\infty(\cF_T)$. But on the other hand, in light of the previous
works (see, e.g., \cite{CHMP} and \cite{Pln}), we see that
technically the domain of $\cE$ should also include the following
set:
 \bea
 \label{cLinf}
 \sL^\infty_T\dfnn\{\xi=\xi_0+zB_T:~\xi_0\in L^\infty(\cF_T),
  ~ z\in \hR^d\}.
 \eea
Here $B$ is the driving Brownian motion. A simple observation of the
Axioms (A3) and (A4) clearly indicates that $\cE$ cannot be defined
simply as  a mapping from $\sL^\infty_T$ to $\sL^\infty_t$. For
example, in general the random variable ${\bf 1}_A\xi$ will not even
be an element of $\sL^\infty_T$(!), thus (A4) will not make sense.

To overcome this difficulty let us now find a larger subset
$\L\subseteq L^0(\cF_T)$ that contains $\sL^\infty_T$ and can serve
as a possible domain of a nonlinear expectation.
First, we observe that such a set must satisfy the following
property in order that Axioms (A1)--(A4) can be well-defined.

\ms

\begin{defn}
Let $\sD(\cF_T)$ denote the totality of all subsets $\L$ in
$L^0(\cF_T)$ satisfying: for all $t\in [0,T]$, the set
 $\L_t \dfnn \L \cap L^0(\cF_t)$ is closed under the multiplication with
 $\cF_t$ indicator functions. That is, if $\xi \in \L_t$ and
 $A \in \cF_t$, then $\b1_A \xi \in \L_t$. \qed
\end{defn}

It is easy to see that $L^\infty(\cF_T) \in \sD(\cF_T)$ and
$\sD(\cF_T)$ is closed under intersections and unions. Thus for any
$S \subset L^0(\cF_T)$, we can define the smallest element in
$\sD(\cF_T)$ that contains $S$ as usual by $\dis \L(S)\dfnn \neg
\neg \underset{\L \in \sD(\cF_T),\, S \subset \L}{\bigcap} \neg \neg
\L$. We are now ready to define the quadratic $\cF$-consistent
nonlinear expectations.

\begin{defn}
\label{fcne}
 An $\cF$-consistent nonlinear expectation
with domain  $\L$ is a pair $(\cE, \L)$, where $\L \in
\sD(\cF_T)$, and $\cE=\{\cE_t\}_{t\ge0}$ is a family of operators
$\cE_t: \L \mapsto \L_t$, $t \in [0,T]$, satisfying Axioms
(A1)--(A4).

Moreover, $\cE$ is called ``translation invariant" if $\L+L^\infty_T \subset \L$ and (\ref{transinv}) holds for any $\xi \in \L$, any $t \in [0,T]$ and any $\eta \in L^\infty(\cF_t)$. \qed
\end{defn}
Again, we shall denote $\cE_t[\cd]=\cE[\cdot|\cF_t]$ as usual, and
we denote $\L=Dom(\cE)$ to be the domain of $\cE$. To simplify
notations, in what follows when we say an $\cF$-consistent
nonlinear expectation $\cE$, we always mean the pair
$(\cE,Dom(\cE))$. Note that a standard $g$-expectation and the
$\cF$-consistent nonlinear expectation studied in \cite{CHMP} and
\cite{Pln} all have domain $\L=L^2(\cF_T)$, and they are
translation invariant if $g$ is independent of $y$. The quadratic
$g$-expectation studied in \cite{MY} is one with domain
$\L=L^\infty(\cF_T)$.

We now turn to the notion of ``quadratic" $\cF$-consistent
nonlinear expectations.
 \begin{defn}
 \label{qne}
An $\cF$-consistent nonlinear expectation $(\cE, Dom(\cE))$ is
called upper (resp. lower) semi-quadratic if there exists a
quadratic $g$-expectation $(\cE^g,Dom(\cE^g))$ with $Dom(\cE^g)
\subseteq Dom(\cE)$ such that for any $t \in [0,T]$ and any $\xi \in
Dom(\cE^g)$, it holds that
 \bea
 \cE[\xi|\cF_t] \leq (\mbox{resp.} \geq) ~\cE^g[\xi|\cF_t], \q~~ \pas
 \eea
Moreover, $\cE$ is called quadratic if there exist two quadratic
$g$-expectations $\cE^{g_1}$ and $\cE^{g_2}$ with
$Dom(\cE^{g_1})\cap Dom(\cE^{g_2}) \subseteq Dom(\cE)$, such that
for any $t \in [0,T]$ and any $\xi \in Dom(\cE^{g_1})\cap
Dom(\cE^{g_2})$, it holds that
 \bea
 \label{gdom*}
  \cE^{g_1}[\xi|\cF_t] \leq \cE[\xi|\cF_t] \leq \cE^{g_2}[\xi|\cF_t], \q~~
  \pas
  \eea
\end{defn}

In what follows, we shall call an $\cF$-consistent nonlinear
expectation as an ``$\cF$-expectation" for simplicity. Note that a
quadratic $g$-expectation $(\cE^g,L^\infty(\cF_T))$ would be a
trivial example of quadratic $\cF$-expectations. The following
example is a little more subtle.

\begin{eg}
\label{cBSDE} {\rm Consider the BSDE (\ref{BSDE}) in which the
generator $g$ is Lipschitz in $y$ and has quadratic growth in $z$.
Furthermore, assume that $g$ is convex in $(t,y,z)$. Then, by a
recent result of Briand-Hu \cite{BH-07}, for any $\xi \in
L^0(\cF_T)$ such that it has exponential moments of all orders (i.e.
$\dis E\big\{e^{\l|\xi|}\Big\}<\infty$, $\fa \l>0$), the BSDE
(\ref{BSDE}) admits a unique solution $(Y,Z)$.
In particular, if we assume further that $g$
satisfies $g|_{z=0}=0$, then it is easy to check the
$g$-expectation $\cE^g(\xi)=Y_0$ defines an
$\cF$-expectation with domain $Dom(\cE^g)\dfnn \Big\{\xi \in \cF_T:
E\Big[e^{\l|\xi|}\Big]<\infty,~\fa \l>0\Big\}$. We should note that
in this case the domain indeed contains the set $\sL^\infty_T$
defined in (\ref{cLinf})!
}\qed
\end{eg}

Since we are only interested in the quadratic $g$-expectations whose
domain contains at least the set $\sL^\infty_T$, we now introduce
the following notion.
\begin{defn}
A quadratic $g$-expectation $\cE^g$ is called ``regular" if
 \beas
\{\xi + zB_\t: \xi \in
 L^\infty(\cF_T),~z \in \hR^d, ~\t \in \cM_{0,T}\}\subseteq Dom(\cE^g).
 \eeas
Correspondingly, a  (semi)-quadratic $\cF$-expectation is called
``regular" if it is dominated by regular quadratic $g$-expectation
in the sense of Definition \ref{qne}.
\end{defn}

Example \ref{cBSDE} shows the existence of the regular quadratic
$g$-expectations. But it is worth pointing out that because of
special form of the set $\sL^\infty_T$, the class of regular
quadratic $g$-expectations is much larger. To see this, let us
consider any quadratic BSDE with $g$ satisfying (H1) and (H2),
 \bea
 \label{qBSDE}
 Y_t=\xi+z B_\t +\int_t^T g(s,Z_s)ds -\int_t^T Z_s dB_s, \q~~ t \in [0,T],
 \eea
where $\xi\in L^\infty(\cF_T)$, $z\in \hR^d$, and $\t \in
\cM_{0,T}$.
Now, if we set $\wt{Y}_t=Y_t-zB_{t\land\t},~\wt{Z}_t=Z_t-z\b1_{\{t
\leq \t \}}$, then (\ref{qBSDE}) becomes
 \bea
 \label{qBSDE1}
 \wt{Y}_t=\xi+\int_t^T g(s,\wt{Z}_s+z\b1_{\{s \leq \t \}})ds -\int_t^T \wt{Z}_s dB_s,
 \q~~ \fa t \in [0,T].
 \eea
Since $\xi\in L^\infty(\cF_T)$, the BSDE (\ref{qBSDE1}) is
uniquely solvable whenever $g$ satisfies (H1) and (H2). In other
words, any $g$ satisfying (H1) and (H2) can generate a regular
$g$-expectation!
%
\begin{rem}
\label{extLinf} {\rm For any generator $g$ satisfying (H1) and
(H2), one can deduce in the similar way as in (\ref{qBSDE}) and
(\ref{qBSDE1}) that
 \bea
 \label{tildeLinf}
 \tilde{\sL}^\infty_T
 \dfnn \Big\{\xi+\int_0^T \z_s dB_s: \xi \in L^\infty(\cF_T), ~ \z
 \in L^\infty_\bF([0,T];\hR^d) \Big\} \subset Dom(\cE^g).
 \eea
Therefore, it follows from Definition \ref{qne} that
 $ \tilde{\sL}^\infty_T
 \subset Dom(\cE^{g_1})\cap Dom(\cE^{g_2})\subset Dom(\cE),
 $
as both $g_1$ and $g_2$ satisfy (H1) and (H2). The set
$\tilde{\sL}^\infty_T$ is very important for the proof of
representation theorem in the last section.
\qed }
\end{rem}

\no{\bf Domination of quadratic $\cF$-expectations.}

\ms

In the theory of nonlinear expectations, especially in the proofs of
decomposition and representation theorems (cf. \cite{CHMP} and
\cite{Pln}), the notion ``domination" for the difference of two
values of $\cF$-expectations plays a central role. To be more
precise, it was assumed that the following property holds for an
$\cF$-expectation $\cE$: for some
$g$-expectation $\cE^g$, it holds for any $X, Y\in L^2(\cF_T)$
that
 \bea
 \label{gdom}
 \cE(X+Y)-\cE(X)\le \cE^g(Y).
 \eea
In the case when $g$ is Lipschitz, this definition of domination
is very natural (especially when $g=g(z)=\m |z|$, $\m>0$).
However, this notion becomes very ill-posed in the quadratic case.
We explain this in the following simple example.
 \begin{eg}
 \label{exdom}
 {\rm Consider the simplest quadratic case: $g=g(z)=\frac{1} {2}|z|^2$,
 and take $\cE=\cE^g$. We show that even such a simple quadratic
 $g$-expectation cannot find a domination in the sense of
 (\ref{gdom}). Indeed, note that
 \beas \cE^g(X+Y)&=&X+Y+\frac{1}{2}
 \int_0^T |Z^1_s|^2ds-\int_0^T Z^1_s dB_s;\\
 \cE^g(X)&=&X+\frac{1}{2}\int_0^T |Z^2_s|^2ds-\int_0^T Z^2_s dB_s.
 \eeas
Denoting $Z=Z^1-Z^2$ we have
$$\cE^g(X+Y)-\cE^g(X)=Y+\frac{1}{2}\int_0^T
(|Z^2_s+Z_s|^2-|Z^2_s|^2)ds-\int_0^T Z_s dB_s.
$$
But in the above the drift $\frac12(|Z^2_s+Z_s|^2-|Z^2_s|^2)\le
|Z_s|^2+\frac12|Z^2_s|^2$ cannot be dominated by any $g$ satisfying
(H1) and (H2).} \qed
\end{eg}

Since finding a general domination rule in the quadratic case is a
formidable task, we are now trying to find a reasonable
replacement that can serve our purpose. It turns out that the
following definition of domination is sufficient for our purpose.
\begin{defn}
\label{Lpdom}
1) A regular quadratic  $\cF$-expectation $\cE$ is said to satisfy
the ``$L^p$-domination" if  for any $K,R
>0$, there exist constants $p=p(K,R)>0$ and $C=C_R>0$ such that for
any two stopping times $0 \leq \t_2 \leq \t_1 \leq T$, any $\xi_i
\in L^\infty_{\t_i}$ with $\|\xi_i \|_\infty \leq K$, $i=1,2$, and
any $z \in \hR^d$ with $|z| \leq R$, it holds for each $t\in[0,T]$
that
 \bea
 \label{pdom}
 \big\| (\cE\{\xi_1\neg+\neg zB_{\t_1}|\cF_t\}\neg-\neg zB_{t
 \land \t_1})\neg-\neg
 (\cE\{\xi_2\neg+\neg zB_{\t_2}|\cF_t\}\neg-\neg zB_{t \land
 \t_2})\big\|_{p}\le  3
 \| \xi_1-\xi_2 \|_{p}+C_R \|\t_1 - \t_2 \|_{p}.
 \eea

2) A regular quadratic  $\cF$-expectation $\cE$ is said to satisfy
the ``$L^\infty$-domination" if for any stopping time $\t \in
\cM_{0,T}$,  any $\xi_i \in L^\infty(\cF_T)$, $i=1,2$, and any $z
\in \hR^d$,  the process $\{\cE\{\xi_i+zB_\t|\cF_t\}-zB_{t \land
\t}, t \in [0,T]\}\in L^\infty_\bF([0,T])$, $i=1,2, $ and
 \bea
 \label{infdom}
 \big\|\cE\{\xi_1+zB_\t|\cF_t\}-\cE\{\xi_2+zB_\t|\cF_t\}
 \big\|_\infty \leq
 \| \xi_1-\xi_2\|_\infty,\q~ \fa t \in [0,T].
  \eea

3) A regular quadratic  $\cF$-expectation $\cE$ is said to satisfy
the ``one-sided $g$-domination" if for any $ K,R
>0$, there are constants $J=J(K,R)>0$
and $\a=\a(K,R)>0$, such that for any stopping time $\t \in
\cM_{0,T}$, $\xi \in L^\infty(\cF_T)$ with $\|\xi \|_\infty \leq K$,
and any $z \in \hR^d$ with $|z| \leq R$, there is a $\g \in BMO$
with $\|\g \|^2_{BMO} \leq J(K,R)$ and a function $g_\a(z)\dfnn
\a(K,R)|z|^2$, $z\in\hR^d$,
such that for any $\eta \in L^\infty(\cF_T)$, it holds that
 \bea
 \label{1sdom}
 \cE[\eta+\xi+zB_\t|\cF_t]-\cE[\xi+zB_\t|\cF_t] \leq
 \cE^{g_\a}_\g[\eta|\cF_t],\q~~\fa t \in [0,T], \q
 P^\g\mbox{-a.s.}
 \eea
Here, $P^\g$ is defined by $dP^\g/dP=\sE(\g \bullet B)_T$, and
$\cE^{g_\a}_\g$ is the $g_\a$-martingale on the probability space
$(\O, \cF, P^\g)$, and with Brownian Motion $B^\g$.
\qed
\end{defn}

 The following theorem more or less justifies the ideas of these ``dominations".

\begin{thm}
\label{justdom} Assume that $g$ is a random field satisfying (H1)
and (H2), and that it satisfies
$g|_{z=0}=0$.
Then the quadratic $g$-expectation $\cE^g$ satisfies both $L^p$
and $L^\infty$-dominations (\ref{pdom}) and (\ref{infdom}).

Furthermore, if $g$ also satisfies that
$\dis\Big|\frac{\pa^2 g}{\pa z^2}\Big| \leq \ell'$ for some
$\ell'>0$, then $\cE^g$ also satisfies the one-sided $g$-domination
(\ref{1sdom}) with $\a(K,R)\equiv \ell'/2$.
\end{thm}

{\it Proof.} 1) We first show that the $L^p$-domination holds. Let
$(Y^i,Z^i)$, $i=1,2$ be the unique solution of BSDE (\ref{qBSDE})
for $\xi_i+zB_{\t_i}$, $i=1,2$, respectively. Define $U_t^i \dfnn
Y_t^i-zB_{t \land \t_i}$, $V_t^i=Z_t^i - z{\bf 1}_{\{t\leq \t_i\}}$,
$\D U_t=U^1-U^2$, and $\D V=V^1-V^2$. Then, in light of
(\ref{qBSDE}) and (\ref{qBSDE1}) one can easily check that
 \bea
 \label{DeltaU1}
 \D U_t=\xi_1-\xi_2+\int_{t \vee \t_2}^{t \vee \t_1}
 g(s,z)ds+\int_t^T \lan\g_s,\D V_s \ran ds-\int_t^T \D V_s
 dB_s,\q \fa t \in [0,T],
 \eea
where $\dis \g_t \dfnn \b1_{\{t \leq \t_1\}}\int_0^1 \frac{\pa
g}{\pa z}\big(t,V^2_t+\th \D V_t+z\big)d\th$. In what follows we
shall denote all the constants depending only on $T$ and $\ell$ in
(H2) by a generic one $C>0$, which may vary from line to line.
Applying Proposition \ref{BMO1} and Corollary 2.2 of \cite{Ko}
we see that both $V^1$ and $V^2$ are BMO with
$$\|V^i\|^2_{BMO}\le C\exp\{C(1+|z|^2)[1+|z|^2+\|\xi_i\|_\infty]\}.
$$
Thus, by definition of $\g$ we have, for any $K, R>0$, with
$\|\xi^1\|_\infty\vee\|\xi^2\|_\infty\le K$ and $|z|\le R$,
 \bea
 \label{gamma}
 \|\g \|^2_{BMO}
 &\leq& C\big[1+|z|^2+\| V^1 \|^2_{BMO}+\| V^2 \|^2_{BMO} \big] \nonumber\\
 &\leq& C(1+|z|^2)+C\exp{\Big\{C(1+|z|^2)\big[\|\xi_1 \|_\infty \vee
 \|\xi_2 \|_\infty +1+|z|^2\big]\Big\}}\\
 &\leq& C(1+R^2)+C\exp{\Big\{C(1+R^2)\big[1+K+R^2\big]\Big\}}\dfnn
 J(K,R).\nonumber
 \eea
Let us now denote $\sE(\g)_s^t\dfnn \frac{\sE(\g\bullet
B)_t}{\sE(\g\bullet B)_s}= \exp\big\{\int_s^t \g_r
dB_r-\frac12\int_s^t|\g_s|^2ds\big\}$, for $0\le s\le t$, and define
a new probability measure $P^\g$ by $dP^\g/dP\dfnn\sE(\g)_0^T$.
Since $\g$ is BMO, applying the Girsanov Theorem we derive from
(\ref{DeltaU1}) that
 \bea
 \label{DeltaU}
 \D U_t=E^\g\Big\{\xi_1-\xi_2+\int_{t \vee \t_2}^{t \vee \t_1}
 g(s,z)ds\Big|\cF_t\Big\} =E\Big\{\big(\xi_1-\xi_2+\int_{t \vee \t_2}^{t
 \vee \t_1} g(s,z)ds\big)\sE(\g)_t^T\Big|\cF_t\Big\},
 \eea
for all $t \in [0,T]$. Since $g$ satisfies (H2), applying the
H\"older inequality we have, for any $p,q>1$ with $1/p+1/q=1$,
 $$|\D U_t|^p\le
 E\Big\{[|\xi_1-\xi_2|+\ell(1+|z|^2)|\t_1-\t_2|]^p|\cF_t\Big\}
 E\Big\{\big[\sE(\g)_t^T\big]^q\big|\cF_t\Big\}^{p/q}.
 $$
Now recall the function $\f_\a$ defined by (\ref{phialpha}). Let
$\a=3$ and $q=q(K,R)>1$ so that $\f_3(q)=J(K,R)$. Applying the
Reversed H\"older Inequality (\ref{rhi}) we obtain, for
$p=p(K,R)=q/(q-1)$,
$$|\D U_t|^p\le 3^p
E\Big\{[|\xi_1-\xi_2|+\ell(1+|z|^2)|\t_1-\t_2|]^p|\cF_t\Big\}.
$$
Taking the expectation, denoting $C_R=3\ell(1+R^2)$, and recalling
the definition of $U$, we have
$$
\big\|(\cE^g[\xi_1+zB_{\t_1}|\cF_t]-zB_{t \land
\t_1})-(\cE^g[\xi_2+zB_{\t_2}|\cF_t] -zB_{t \land \t_2})\big\|_p
\leq 3\|\xi_1-\xi_2\|_p+C_R\|\t_1-\t_2 \|_p,
$$
for all $t \in [0,T]$, proving (\ref{pdom}).

2) The proof of ``$L^\infty$-domination" (\ref{infdom}) is similar
but much easier. Again we let $(Y^i,Z^i)$ be the solution of
(\ref{qBSDE}) for $\xi_i+zB_\t$, $i=1,2$, respectively.
Denote $\D Y=Y^1-Y^2$ and $\D Z=Z^1-Z^2$, we have
 \beas
 \D Y_t=\D \xi+\int_t^T \lan \g_s,\D Z_s\ran ds-\int_t^T \D Z_s
 dB_s,\q~~ \fa t \in [0,T],
 \eeas
where $\dis \g_t \dfnn \int_0^1 \frac{\pa g}{\pa z}\big(t,\l
Z^1_t+(1-\l)Z^2_t\big)d\l \in BMO$. Applying Girsanov's Theorem
again we obtain that, under some equivalent probability measure
$P^\g$, it holds that
 \beas \D Y_t=E^\g [\D \xi |\cF_t],\q~~ \fa t
 \in [0,T], \q \pas
 \eeas
The estimate (\ref{infdom}) then follows immediately.

\ms

3) We now prove the one-sided $g$-domination (\ref{1sdom}). This
times we let $(Y^1,Z^1)$ and $(Y^2,Z^2)$ be the solutions of BSDE
(\ref{qBSDE}) with terminal conditions $\eta+\xi+zB_\t$ and
$\xi+zB_\t$, respectively. Then (\ref{qBSDE1}) implies that, for all
$t \in [0,T]$,
 \beas
 \D Y_t &=& \D \wt{Y}_t= \eta+\int_t^T \Big(g\big(s,\wt{Z}^1_s+z\b1_{\{s \leq \t\}}\big)
  - g\big(s,\wt{Z}^2_s+z\b1_{\{s \leq \t\}}\big)\Big)ds - \int_t^T \D \wt{Z}_s dB_s \\
 &=& \eta+\int_t^T \lan \int_0^1 \frac{\pa g}{\pa z}\big(s,\l\D \wt{Z}_s+\wt{Z}^2_s
 +z\b1_{\{s \leq \t\}}\big)d\l, \D \wt{Z}_s \ran ds- \int_t^T \D
 \wt{Z}_sdB_s,
 \eeas
where $\wt{Y}^i_t \dfnn Y^i_t-zB_{t\land\t}$ and $\wt{Z}^i_t\dfnn
Z^i_t-z\b1_{\{t \leq \t \}}$, $i=1,2$. Since $\wt{Z}^i \in$ BMO,
$i=1,2$, thanks to  Proposition \ref{BMO1}, it is easy to check
that $\g_\cd \dfnn \frac{\pa g}{\pa z}(\cd,Z^2_\cd)\in$ BMO as
well, and the estimate (\ref{gamma}) remains true.
It is worth noting that $\g$ is independent of $\eta$ since $Z^2$
is so. By Girsanov's Theorem,
 \beas
 \D Y_t &=& \eta+\int_t^T \lan \int_0^1 \Big(\frac{\pa g}{\pa z}\big(s,\l\D \wt{Z}_s+
 \wt{Z}^2_s+z\b1_{\{s \leq \t\}}\big)-\frac{\pa g}{\pa z}\big(s,\wt{Z}^2_s+z\b1_{\{s
 \leq \t\}}\big)\Big)d\l, \D \wt{Z}_s \ran ds \\
 &&-\int_t^T \D \wt{Z}_s dB^\g_s,\q~~ \fa t \in [0,T],
 \eeas
where $P^\g$ is the equivalent probability measure as before. Now
the extra assumption on the boundedness of $\frac{\pa^2 g}{{\pa
z}^2}$ concludes that, with $\a(K,R)\equiv \ell'/2$,
 $$
 \bigg|\Big\lan \int_0^1 \Big(\frac{\pa g}{\pa z}\big(s,\l\D \wt{Z}_s+\wt{Z}^2_s+z\b1_{\{s
 \leq \t\}}\big)-\frac{\pa g}{\pa z}\big(s,\wt{Z}^2_s+z\b1_{\{s \leq \t\}}\big)\Big)d\l,
 \D \wt{Z}_s \Big\ran\bigg|
 \leq \a(K,R) \big|\D \wt{Z}_s \big|^2.
 $$
The Comparison Theorem of quadratic BSDE (cf. \cite[Theorem
2.6]{Ko}) then leads to that
 $$
 \cE^g[\eta+\xi+zB_\t|\cF_t]-\cE^g[\xi+zB_\t|\cF_t] \leq
 \cE^{g_\a}_\g[\eta|\cF_t],\q~~\fa t \in [0,T],
 $$
proving (\ref{1sdom}), whence the theorem. \qed

\section{Properties of Quadratic $\cF$-expectations}
\setcounter{equation}{0}

In this section, we assume that $\cE$ is a translation invariant
semi-quadratic $\cF$-expectation dominated by a quadratic
$g$-expectation $\cE^g$ with $g$ satisfying (H1) and (H2). Clearly
$\cE$ is regular. We also assume that $\cE$ satisfies both the
$L^p$-domination (\ref{pdom}) and the $L^\infty$-domination
(\ref{infdom}).

 We first give a path regularity result for $\cE$-martingales, which is very useful
in our future discussion.

\begin{prop} \label{cem}
For any $\t \in \cM_{0,T}$, $\xi \in L^\infty(\cF_\t)$, and  $z \in
\hR^d$, the process $\cE[\xi+ zB_\t|\cF_t]$, $t \in[0,T]$ admits a
\cad modification.
\end{prop}

{\it Proof.} We first assume that $\cE$ is an upper semi-quadratic
$\cF$-expectation first. By the $L^\infty$-domination, $X_\cdot
\dfnn \cE[\xi+zB_\t |\cF_\cdot]-zB_{\cdot \land \t}\in
L^\infty_\bF([0,T])$, which implies that $|X_t|\le \|X\|_\infty$,
$\pas$ for any $t \in [0,T]$ except a null set $\cT$. We may assume
that there is a dense set $\cD$ of $[0,T]\backslash \cT$ such that
$|X_t|\le \|X\|_\infty$, $\fa t \in \cD$, $\pas$ Now we define a new
generator
 \bea
 \label{hatg}
 \hat{g}(t,\o,\z)\dfnn g(t,\o,\z+\b1_{\{t \leq \t\}}z)-g(t,\o,\b1_{\{t \leq
 \t\}}z),\q~~\fa (t,\o,\z) \in [0,T] \times \O \times \hR^d.
 \eea
%
%
For any $0 \leq s \leq t \leq T$ and any $\eta \in L^\infty(\cF_t)$,
it is easy to check that $P$-a.s.
 \bea
 \label{get}
 \cE^g[\eta+zB_{t \land \t}|\cF_s]-zB_{s \land \t}+\int_0^s g(r,\b1_{\{r \leq \t\}}z)dr
 =\cE^{\hat{g}}[\eta+\int_0^t g(r,\b1_{\{r \leq \t\}}z)dr|\cF_s].
 \eea
In particular, by the definition and the properties of upper
semi-quadratic $\cF$-expectation, letting $\eta=X_t$ in (\ref{get})
shows that $P$-a.s.
 \beas
 \cE[\xi+zB_\t|\cF_s] &=&\cE\big[\cE[\xi+zB_\t|\cF_t]\big|\cF_s\big] =
 \cE[X_t+zB_{t \land \t}|\cF_s] \leq \cE^g[X_t+zB_{t \land \t}|\cF_s] \\
 &=& \cE^{\hat{g}}\Big\{X_t+\int_0^t g(r,\b1_{\{r \leq \t\}}z)dr\Big|\cF_s\Big\}+zB_{s \land \t}-\int_0^s
 g(r,\b1_{\{r \leq \t\}}z)dr.
 \eeas
In other words, the process $t\mapsto X_t+\int_0^t g(r,\b1_{\{r \leq
\t\}}z)dr$ is in fact a $\hat{g}$-submartingale. Thus by Theorem
\ref{mayaothm} we can define a \cad process
 \beas
 Y_t \dfnn \underset{r \searrow t,\; r \in \cD}{\lim} X_r, \q \fa t \in [0,T) \q
 \mbox{and}\q Y_T \dfnn X_T= \xi.
 \eeas
Clearly,
$Y \in \hD^\infty_\bF([0,T])$. Moreover, the constant-preserving
property of $\cE$ and ``Zero-One Law" imply that \bea \label{tau_t}
\cE[\xi'|\cF_t] \in \cF_{t \land \t}, \qq \fa \xi' \in \L_\t, \q~~
\fa t \in [0,T]. \eea To see this, one needs only note that for any
$s \in [0,t),$ \beas &&\b1_{\{t \land \t \le s\}}
\cE[\xi'|\cF_t]=\b1_{\{\t \le s\}} \cE[\xi'|\cF_t] = \cE[\b1_{\{\t
\le s\}}\xi'|\cF_t]=\b1_{\{\t \le s\}}\xi' \in \cF_s. \eeas Thus
$X_t \in \cF_{t \land \t},\fa t \in [0,T]$, so is $Y$ by the
right-continuity of the filtration $\bF$. Now, for any $t \in[0,T)$
and $r\in (t,T]\cap \cD$, we write
 \beas
 X_t-Y_t=\cE[\xi+zB_\t |\cF_t]-zB_{t \land \t}-Y_t=\cE[X_r+zB_{r
 \land \t}|\cF_t]-\cE[Y_t+zB_{t \land \t}|\cF_t].
 \eeas
Then applying (\ref{pdom}) with
$K=\|X \|_\infty$ and $R=|z|$ we can find a $p=p(K,R)$ such that
 \beas
 \|X_t-Y_t \|_p\leq 3\| X_r-Y_t \|_p + C_R\| r \land \t - t \land \t
 \|_p\leq 3\| X_r-Y_t \|_p + C_R(r - t).
 \eeas
Letting $r\searrow t$ in the above, the Bounded Convergence Theorem
then implies that $X_t=Y_t$, $\pas$ To wit, the process $Y_t+zB_{t
\land \t}$, $t \in [0,T]$ is a \cad modification of $\cE[\xi+zB_\t
|\cF_t]$, $t \in [0,T]$.

The case when $\cE$ is lower semi-quadratic can be argued similarly.
The proof is complete. \qed



Next, we prove the ``optional sampling theorem" for the quadratic
$\cF$-expectation. To begin with, we recall that the nonlinear
conditional expectation $\cE[\cd|\cF_\si]$ is defined as follows. If
$\xi\in Dom(\cE)$, denote $Y_t\dfnn \cE[\xi|\cF_t]$, $t \in [0,T]$,
then for any $\si\in \cM_{0,T}$, we define
 \bea
 \label{taucond}
 \cE[\xi|\cF_\si] \dfnn Y_\si,\qq \pas
 \eea
The following properties of $\cE[\cd|\cF_\si]$ are important.
 \begin{prop}
 \label{psfe}
For any $\t,\si \in \cM_{0,T}$, $\xi,\eta \in L^\infty(\cF_\t)$, and
$z \in \hR^d$, it holds that
 \bi
 \item[(i)]
$\cE[\xi+zB_\t|\cF_\si] \leq \cE[\eta+zB_\t|\cF_\si],~\pas$, if $\xi
\leq \eta~\pas$;

 \item[(ii)] $\cE[\xi+zB_\t|\cF_\t] = \xi+zB_\t,~ \pas$;

 \item[(iii)] $\b1_A
\cE[\xi+zB_\t|\cF_\si]=\b1_A\cE[\b1_A\xi+zB_\t|\cF_\si],~\pas$, $\fa
A \in \cF_{\t
 \land \si}$;

\item[(iv)] If further $\eta \in L^\infty(\cF_{\t
 \land \si})$, the
following ``translation invariance" property holds:
 \beas
 \cE[\xi+zB_\t+\eta|\cF_\si]=\cE[\xi+zB_\t|\cF_\si]+\eta,  \q~~ \pas
 \eeas
\ei
\end{prop}

{\it Proof.} (i) is a direct consequence of the monotonicity of
$\cE$ and Proposition \ref{cem}.

To see (ii), we first assume that $\t$ takes values in a finite set:
$0\le t_1 <\cdots< t_n\le T$. Actually, for any $\xi' \in \L_\t$,
the constant-preserving of $\cE$ and ``Zero-One Law" imply that
 \beas
\cE[\xi'|\cF_\t]=\sum_{j=1}^n \b1_{\{\t = t_j\}}\cE[\xi'|\cF_{t_j}]
=\sum_{j=1}^n\cE[\b1_{\{\t = t_j\}}\xi'|\cF_{t_j}]=\sum_{j=1}^n
\b1_{\{\t = t_j\}}\xi'=\xi',\q~~\pas
\eeas For general stopping time
$\t$, we first choose a sequence of finite valued stopping times
$\{\t_n\}$ such that
$\t_n \searrow \t$, $\pas$ Since for each
$n$ it holds that
 \beas
 \cE[\xi+zB_\t|\cF_{\t_n}]=\xi+zB_\t, \q~~\pas, \q n=1,2,\cds,
 \eeas
letting $n \to \infty$ and applying Proposition \ref{cem} we
obtain that
$\cE[\xi+zB_\t|\cF_\t]=\xi+zB_\t$, $\pas$, proving (ii).

We now prove (iii). Again, we assume first that $\si$ takes finite
values in $0\leq t_1 < \cdots < t_n \leq T$. For any $A \in \cF_{\t
 \land \si}$, let $A_j= A \cap\{\si = t_j\}\in \cF_{t_j}$, $1 \leq j \leq
n$. Then it holds $P$-a.s. that
 \beas
 \b1_A\cE[\b1_A\xi+zB_\t|\cF_\si]&=&\sum_{j=1}^n \b1_{A_j}\cE[\b1_A\xi+zB_\t|\cF_{t_j}]
 =\sum_{j=1}^n \cE[\b1_{A_j}\xi+\b1_{A_j}zB_\t|\cF_{t_j}]\\
 &=&\sum_{j=1}^n
 \b1_{A_j}\cE[\xi+zB_\t|\cF_{t_j}]=\b1_A\cE[\xi+zB_\t|\cF_\si].
 \eeas
For general stopping time $\si$, we again approximate $\si$ from
above by a sequence of finite-valued stopping times
$\{\si_n\}_{n\ge 0}$.
Then for any $A \in \cF_{\t
 \land \si} \subset \cF_{\t
 \land \si_n}$, $\fa n \in \hN$, we have
 \beas
 \b1_A\cE[\xi+zB_\t|\cF_{\si_n}]=\b1_A\cE[\b1_A\xi+zB_\t|\cF_{\si_n}],
 \q~~ \pas, \q \fa n \in \hN.
 \eeas
Letting $n\to\infty$ and applying Proposition \ref{cem} again we
can prove (iii).

(iv) The proof is quite similar, thus we shall only consider the
case where $\si$ takes values in a finite set $0\leq t_1 < \cdots
 < t_n \leq T$. In this case we have
 \beas
 \cE[\xi+zB_\t+\eta|\cF_\si]&=&\sum_{j=1}^n \b1_{\{\si = t_j\}}\cE[\xi+zB_\t+\eta|\cF_{t_j}]\\
 &=&\sum_{j=1}^n \cE[\b1_{\{\si = t_j\}}(\xi+zB_\t)+\b1_{\{\si = t_j\}}\eta|\cF_{t_j}] \\
 &=&\sum_{j=1}^n \Big\{\cE[\b1_{\{\si = t_j\}}(\xi+zB_\t)|\cF_{t_j}] +\b1_{\{\si = t_j\}}
 \eta \Big\}\\
 &=&\sum_{j=1}^n \b1_{\{\si = t_j\}}\cE[\xi+zB_\t|\cF_{t_j}] +\sum_{j=1}^n \b1_{\{\si = t_j\}}
 \eta=\cE[\xi+zB_\t|\cF_\si] + \eta.
 \eeas
The third equality is due to the ``translation invariance" of
$\cE$ and $\b1_{\{\si = t_j\}}\eta \in L^\infty(\cF_{t_j})$. The
rest of the proof can be carried out in a similar way as other
cases, we leave it to the interested reader.
The proof is complete. \qed

We now prove an important property of $\cE\{\cd|\cF_t\}$, which we
shall refer to as the  ``{\it Optional Sampling Theorem}" in the
future.
 \begin{thm}
 \label{opts}
 For any $X \in L^\infty_\bF([0,T])$ and $z \in \hR^d$ such that $t\mapsto X_t +z B_t$
 is a right-continuous $\cE$-submartingale (resp. $\cE$-supermartingale or $\cE$-martingale).
 Then for any stopping times $\t,\si \in [0,T]$, it holds that
 \beas
 \cE[X_\t+zB_\t|\cF_\si]\geq
 ~(\mbox{resp. }
 \leq \mbox{ or } = )~X_{\t \land \si}+zB_{\t \land \si},\qq \pas
 \eeas
\end{thm}

{\it Proof.} We shall consider only the $\cE$-submartingale
case, as the other cases can be deduced easily by
standard argument. To begin with, we assume that $\si \equiv t \in
[0,T]$ and assume that $\t$ takes finite values in $0\leq t_1 <
\cdots < t_N \leq T$. Note that if $t \geq t_N$, then $X_\t+zB_\t
\in \cF_t$ and $\t \land t=\t$, thus
 \beas
 \cE[X_\t+zB_\t|\cF_t]=X_\t+zB_\t =X_{\t \land t}+zB_{\t \land t}, \qq \pas,
 \eeas
thanks to the constant preserving property of $\cE$. We can then
argue inductively to show that the statement holds for $t\ge t_m$,
for all $1\le m\le N$. In fact, assume that for $m \in \{2,\cdots
N\}$
 \bea
 \label{tgetm}
 \cE[X_\t+zB_\t|\cF_t] \geq X_{\t \land t}+zB_{\t \land t}, \q~ \pas\q~~ \fa t \geq
 t_m.
 \label{opts1}
 \eea
Then, again using the translability and the ``zero-one" law, one
shows that for any $t \in [t_{m-1},t_m)$, it holds $P$-a.s. that
 \beas
 \cE[X_\t+zB_\t|\cF_t]&=&\cE\big[\cE[X_\t+zB_\t|\cF_{t_m}]\big|\cF_t\big] \geq
 \cE[X_{\t \land t_m}+zB_{\t \land t_m}|\cF_t] \\
 &=&\cE[\b1_{\{\t \leq t_{m-1}\}}(X_{\t \land t}+zB_{\t \land t})+\b1_{\{\t \geq t_m\}}
 (X_{t_m}+zB_{t_m})|\cF_t\big]\\
 &=&\b1_{\{\t \leq t_{m-1}\}}(X_{\t \land t}+zB_{\t \land t})+\b1_{\{\t \geq t_m\}}
 \cE[X_{t_m}+zB_{t_m}|\cF_t\big]\\
 &\geq &\b1_{\{\t \leq t_{m-1}\}}(X_{\t \land t}+zB_{\t \land t})+\b1_{\{\t \geq t_m\}}
 (X_t+zB_t)\\
 &=&X_{\t \land t}+zB_{\t \land t}.
 \eeas
Namely
(\ref{opts1}) also holds for any $t \geq t_{m-1}$. This completes the inductive step.
Thus (\ref{opts1}) holds for all finite-valued stopping times.

Now let $\t$ be a general stopping time, we still choose $\{\t_n\}$
to be a sequence of finite-valued stopping times such that
$\t_n\searrow \t$, $\pas$ Then (\ref{opts1}) holds for all $\t_n$'s.
Now let $K=\|X \|_\infty$, $R = |z|$, and $p=p(K,R)$. Applying the
$L^p$-domination (\ref{pdom}) for $\cE$ we see that for any $n \in
\hN$, \bea \label{opts3}
& & \|\cE[X_{\t_n}+zB_{\t_n}|\cF_t]-\cE[X_\t+zB_\t|\cF_t] \|_p\nonumber \\
&\leq & \|(\cE[X_{\t_n}+zB_{\t_n}|\cF_t]-zB_{\t_n \land t})-(\cE[X_\t+zB_\t|\cF_t]
-zB_{\t \land t}) \|_p+R \|B_{\t_n \land t}-B_{\t \land t} \|_p \nonumber \\
&\leq & 3\| X_{\t_n}-X_\t \norm_p+C_R \norm \t_n - \t\|_p+ R \|
B_{\t_n \land t}-B_{\t \land t} \|_p. \eea Since $X$ is a bounded
\cad process, we can then apply the Bounded Convergence Theorem to
conclude that the first and second terms on the right hand side of
(\ref{opts3}) tend to $0$, as $n\to \infty$.  Furthermore, applying
the Burkholder-Davis-Gundy inequality and Bounded Convergence
Theorem, we conclude that the last term on the right hand side of
(\ref{opts3}) also goes to $0$. Thus, possibly along a subsequence,
we see that for any $t \in [0,T]$ \beas \cE[X_\t+zB_\t|\cF_t]=\neg
\underset{n \to \infty}{\lim}\cE[X_{\t_n}+zB_{\t_n}|\cF_t] \ge
\underset{n \to \infty}{\lim} \big(X_{\t_n \land t}+zB_{\t_n \land
t}\big)=X_{\t \land t}+zB_{\t \land t}, \q \pas \eeas
Thus we obtain (\ref{opts1}) again.

Finally, let us consider the case when $\si$ is also a general
stopping time. Following the previous argument, with the help of
Proposition \ref{cem}, we have, $P$-a.s.
 \beas
 \cE[X_\t+zB_\t|\cF_t] \geq X_{\t \land t}+zB_{\t \land t}, \qq \fa
 t \in [0,T].
 \eeas
Consequently, we obtain that $\cE[X_\t+zB_\t|\cF_\si] \geq X_{\t
\land \si}+zB_{\t \land \si}$, $\pas$, proving the theorem. \qed

To end this section we consider a special BSDE involving the
quadratic $\cF$-expectation $\cE$, which will be very useful in
the rest of the paper:
 \bea
 \label{BSDEcE}
 Y_t+zB_t+\int_0^t f(s,Y_s)ds= \cE\Big\{\xi+zB_T+\int_0^T
 f(s,Y_s)ds\Big|\cF_t\Big\}, \qq \fa t \in [0,T],
 \eea
where $f: [0,T]\times \O \times \hR \to \hR$ is a measurable
function such that it satisfies the following assumption:
 \begin{itemize}
\item[{\bf (H3)}] The function $f$ is uniformly Lipschitz in $y$
with Lipschitz constant $\k>0$, uniform in $(t,\o)$, such that
$\int_0^T |f(t,\cd,0)|dt \in L^\infty(\cF_T)$.
 \end{itemize}

We have the following existence and uniqueness result for the BSDE
(\ref{BSDEcE}).
\begin{prop} \label{qfee}
Assume (H3). Then for any $\xi \in L^\infty(\cF_T)$ and any $z \in
\hR^d$, the BSDE \(\ref{BSDEcE}\) admits a unique solution in
$\hD^\infty_\bF([0,T])$.
\end{prop}


{\it Proof.} We first consider the case when $T\le 1/2\k$, where
$\k$ is the Lipschitz constant of $f$ in (H3). For any $Y \in
\hD^\infty_\bF([0,T])$, and $t \in [0,T]$, using (H3) we have
$$ \Big\|\int_0^t f(s,Y_s) ds \Big\|_\infty \leq \Big\| \int_0^T
|f(s,0)|ds \Big\|_\infty +\k t\|Y \|_\infty < \infty.
$$
In particular, we have $\xi+\int_0^T f(s,Y_s)ds \in L^\infty(\cF_T)$
so that $\cE\{\xi+zB_T+\int_0^T f(s,Y_s)ds|\cF_t\}$ is well-defined,
and we can define a mapping $\F: \hD^\infty_\bF([0,T]) \mapsto
\hD^\infty_\bF([0,T])$ by:
 \bea
 \label{Phi}
 \F_t(Y)\dfnn \cE\Big\{\xi+zB_T+\int_0^T f(s,Y_s)ds\Big|\cF_t\Big\}-zB_t
 -\int_0^t f(s,Y_s)ds, \qq t\in[0,T].
 \eea

We claim that $\F$ is a contraction. Indeed, since $\cE$ satisfies
the $L^\infty$-domination,
for any $Y,\hat{Y} \in \hD^\infty_\bF([0,T])$, (\ref{infdom})
implies that for any $t \in [0,T]$, it holds $P$-a.s. that
 \bea
 \label{qfee1}
 |\F_t(Y)-\F_t(\hat{Y})|\neg &\neg\neg=\neg\neg& \neg
 \Big|\cE[\xi+zB_T+\int_t^T f(s,Y_s)ds|\cF_t]
 -\cE[\xi+zB_T+\int_t^T f(s,\hat{Y}_s)ds|\cF_t]\Big|\qq \nonumber \\
\neg&\neg\neg\leq\neg\neg & \neg \Big\|\neg \int_t^T \neg\neg
\big(f(s,Y_s)\neg- \neg f(s,\hat{Y}_s)\big)ds \Big\|_\infty \leq \k
(T-t) \| Y-\hat{Y} \|_\infty  \leq \frac{1}{2} \| Y-\hat{Y}
\|_\infty.
 \eea
Since the process  $t\mapsto \F_t(Y)$ is \cad\neg, thanks to
Proposition \ref{cem}, we conclude that
$\|\F(Y)-\F(\hat{Y}) \|_\infty \leq \frac{1}{2} \|Y-\hat{Y}
\|_\infty$. Thus $\F$ is a contraction, and the lemma holds in
this case.

The general case can now be argued using a standard ``patching-up"
method. Namely we take a partition of $[0,T]$: $0=t_0<t_1<\cds
<t_N=T$, such that $\max |t_n-t_{n-1}|<1/2\k$. We first solve the
BSDE (\ref{BSDEcE}) on $[t_{N-1}, t_N]$ to get a solution $Y^N$.
We then solve (\ref{BSDEcE}) on $[t_{N-2},,t_{N-1}]$ to get
$Y^{N-1}$, satisfying the terminal condition
$Y^{N-1}_{t_{N-1}}=Y^N_{t_{N-1}}$, and so on,
thanks to the result proved in the first part. Denoting the
solution on $[t_{n-1}, t_n]$ by $Y^n$,
we can then define a new process by
 $Y_t \dfnn Y^n_t$, $ t \in [t_{n-1},t_n]$, $ n=1,\cds,N$, and
prove that $Y$ solves (\ref{BSDEcE}) over $[0,T]$ by induction.

To see this, we first note that $Y\in \hD^\infty_\bF([0,T])$. Now
assuming that $Y$ solves (\ref{BSDEcE}) on $[t_n,T]$, we show that
it solves (\ref{BSDEcE}) on $[t_{n-1},T]$ as well. Indeed,
for any $t\in[t_{n-1}, t_n]$, we have
%
 \beas
  Y_t+zB_t+\int_0^t f(s,Y_s)ds &=& Y^{n}_t+zB_t+\int_{t_{n-1}}^t f(s,Y^{n}_s)ds +
 \int_0^{t_{n-1}} f(s,Y_s)ds \\
 &=& \cE\Big\{Y^n_{t_n}+zB_{t_n} +\int_{t_{n-1}}^{t_n} f(s,Y^{n}_s)ds\Big|\cF_t\Big\}+
 \int_0^{t_{n-1}} f(s,Y_s)ds \\
 &=& \cE\Big\{Y_{t_n}+zB_{t_n}+\int_0^{t_n} f(s,Y_s)ds \Big|\cF_t\Big\} \\
 &=& \cE\Big\{\cE\Big\{\xi+zB_T+\int_0^T f(s,Y_s)ds \Big|\cF_{t_n}\Big\}\Big|\cF_t\Big\}
  \\
 &=& \cE\Big\{\xi+zB_T+\int_0^T f(s,Y_s)ds \Big|\cF_t\Big\}.
 \eeas
In the above the second equality is due to the fact that $Y^n$
solves (\ref{BSDEcE}) on $[t_{n-1},t_n]$; the third equality is due
to the ``translation invariance" of $\cE\{\cd|\cF_t\}$; the fourth
equality is because of the inductional hypothesis that $Y$ solves
(\ref{BSDEcE}) on $[t_n,T]$; and the last equality is the
``time-consistence" property of $\cE\{\cd|\cF_t\}$. This shows that
$Y$ solves (\ref{BSDEcE}) on $[t_{n-1},T]$, whence the existence.

The uniqueness can be argued in a similar way.
First note that the BSDE (\ref{BSDEcE}) can be written in a
``local" form: for $n=1,2,\cds, N$,
 \bea
 \label{localYn}
 Y_t+zB_t
 =\cE\Big\{Y_{t_n}+zB_{t_{n}}+\int_t^{t_{n}} f(s,Y_s)ds \Big|\cF_t\Big\}, \qq t\in
 [t_{n-1},t_n],
 \eea
thanks to the translation invariance property of $\cE\{\cd|\cF_t\}$.
Assume that $\hat{Y}\in \hD^\infty_\bF([0,T])$ is another solution
of (\ref{BSDEcE}). Then it must satisfy (\ref{localYn}) on
$[t_{N-1}, T]$. The fixed point argument in the first part then
shows that $Y=\hat Y$ in $\hD^\infty_\bF([t_{N-1},T])$,
thus $Y_{t_{N-1}}=\hat Y_{t_{N-1}}$, $\pas$ We can repeat the same
argument for $[t_{N-2}, t_{N-1}]$, and so on to conclude after
finitely many steps that
$Y$ and $\hat{Y}$ are indistinguishable over the whole interval
$[0,T]$. The proof is now complete. \qed

\section{Doob-Meyer Decomposition of Quadratic $\cF$-Martingales}
\setcounter{equation}{0}

In this section we prove a Doob-Meyer type decomposition theorem for
quadratic $\cF$-martingales. We shall assume that $\cE$ is a
translation invariant quadratic $\cF$-expectation dominated by two
quadratic $g$-expectations $\cE^{g_1}$ and $\cE^{g_2}$ from below
and above, and both $g_1$ and $g_2$ satisfies (H1) and (H2) with the
same $\ell>0$. We also assume that $\cE$ satisfies both the
$L^p$-domination (\ref{pdom}) and the $L^\infty$-domination
(\ref{infdom}).

The following proposition will play an essential role in the rest of
this paper.
 \begin{prop}
 \label{cfm}
For any $\t \in \cM_{0,T}$, $\xi \in L^\infty(\cF_\t)$, and $z \in
\hR^d$, denote $Y_t \dfnn \cE[\xi+zB_\t|\cF_t]$, $t \in [0,T]$. Then
there exists a unique pair $\dis (h,Z) \in L^1_\bF([0,T]) \times
\cH^2_\bF([0,T]; \hR^d)$ such that \bea
 \label{hdom}
  -\ell(|Z_t|+|Z_t|^2) \leq g_1(t,Z_t) \leq h_t \leq
 g_2(t,Z_t)\leq \ell(|Z_t|+|Z_t|^2), \qq  dt \times dP\mbox{-a.s.},
 \eea
 and $(Y,Z)$ satisfies the BSDE:
\bea
 \label{cfm1}
Y_t=Y_T+\int_t^T h_s ds -\int_t^T Z_s dB_s, \qq \fa t \in [0,T].
 \eea
Moreover, if we assume that $\cE$ also satisfies the
one-sided $g$-domination (\ref{1sdom}), with $K \ge \| \xi
\|_\infty$, $R \ge |z|$, $\a=\a(K,R)$, $J=J(K,R)$ and $\| \g
\|^2_{BMO}\le J$, then for any $\eta \in L^\infty(\cF_\t)$, the pair
$(\hat{h},\hat{Z})$ corresponding to the process
$\cE\{\eta+zB_\t|\cF_t\}$, $t \in [0,T]$, satisfies
\bea \hat{h}_t-h_t \leq \a|\hat{Z}_t-Z_t|^2+\lan \g_t,\hat{Z}_t-Z_t
\ran, \qq  \dtp \eea
\end{prop}

{\it Proof.} For each $z\in\hR^d$, define a process $\wt{Y}_t \dfnn
Y_t-zB_{t \land \t}$, $t\in [0,T]$ and a new generator
 \beas
 g_i^z(t,\o,\z) \dfnn g_i(t,\o,\z+\b1_{\{t \leq \t\}}z), \q~~ \fa
 (t,\o,\z)\in [0,T] \times \O \times \hR^d,\q \q i=1,2.
 \eeas
By the definition of the $L^\infty$-domination (see Definition
\ref{Lpdom}-(2)) and the fact (\ref{tau_t}) we see that $\wt{Y}
\in L^\infty_\bF([0,T])$ and $\wt{Y}_t \in \cF_{t \land \t}$, $\fa
t \in [0,T]$.
It is easy to check that for $0 \leq s \leq t \leq T$ and any
$\eta \in L^\infty(\cF_t)$, \beas \cE^{g_i}[\eta+zB_{t \land
\t}|\cF_s]=\cE^{g_i^z}[\eta|\cF_s]+zB_{s \land \t}, \qq \pas \q~~
i=1,2. \eeas Thus the upper domination of $\cE$ by $\cE^{g_1}$ and
the time-consistency of $\cE$ imply that, $P$-a.s.,
 \beas
 \cE^{g_1^z}[\wt{Y}_t|\cF_s]&=& \cE^{g_1}[\wt{Y}_t+zB_{t \land
 \t}|\cF_s]-zB_{s \land \t}
 =\cE^{g_1}\big[ \cE[\xi+zB_\t|\cF_t] \big|\cF_s\big]-zB_{s \land \t}\\
 &\leq& \cE\big[ \cE[\xi+zB_\t|\cF_t] \big|\cF_s\big]-zB_{s \land
 \t}=\cE[\xi+zB_\t|\cF_s]-zB_{s \land \t}=\wt{Y}_s.
 \eeas
Namely, $\wt{Y}$ is both a $g^z_1$-supermartingale and a
$g^z_2$-submartingale. Applying Theorem \ref{mayaothm} we
obtain two increasing processes $A^1$ and $A^2$ (we may assume both
are \cad and null at $0$) and two processes $\wt{Z}^1$, $\wt{Z}^2 \in
\cH^2_\cF(\hR^d)$,  such that
 \beas
 \wt{Y}_t=\wt{Y}_T+\int_t^T
 g^z_i(s,\wt{Z}^i_s)ds +(-1)^i( A^i_t-A^i_T)-\int_t^T
 \wt{Z}^i_s dB_s, \q~~ t \in [0,T], \q i=1,2.
 \eeas
Letting $Z^i_t = \wt{Z}^i_t+\b1_{\{t \leq \t\}}z$ we have, for
$i=1,2$,
 \bea\label{2bsde}
 Y_t=Y_T+\int_t^T g_i(s,Z^i_s)ds + (-1)^i(
 A^i_t-A^i_T)-\int_t^T Z^i_s dB_s, \q~~ \fa t \in [0,T].
 \eea
By comparing the martingale parts and bounded variation parts of two
BSDEs in (\ref{2bsde}), one has:
 \beas Z^1_t\equiv Z^2_t, \q
 \mbox{and} \q -g_1(t,Z^1_t)dt-dA^1_t \equiv -g_2(t,Z^2_t)dt+dA^2_t, \q t\in[0,T],\q\pas
 \eeas
Consequently, we have that $dA^1_t+dA^2_t \equiv
\big(g_2(t,Z^1_t)-g_1(t,Z^1_t)\big)dt$, which implies that both
$A^1$ and $A^2$ are absolutely continuous and $dA^i_t=a^i_tdt$ with
$a^i_t \ge 0$, $i=1,2$. The conclusion follows by setting $Z_t \dfnn
Z^1_t$ and $h_t \dfnn g_1(t,Z_t)+a^1_t$.

Moreover, if $\cE$ also satisfies the one-sided $g$-domination
(\ref{1sdom}), then for any $\eta \in L^\infty(\cF_\t)$, we can set
$\hat{Y}_t \dfnn \cE[\eta+zB_\t|\cF_t]$, $\fa t \in [0,T]$ and let
$(\hat{h},\hat{Z})$ be the corresponding pair. Applying the
$L^\infty$-domination (\ref{infdom}) for $\cE$, we see that
$\hat{Y}-Y \in L^\infty_\bF([0,T])$ under $P$, whence under $P^\g$.
In fact, $\hat{Y}-Y$ is a $g_\a$ submartingale under $P^\g:$ for $0
\leq s \leq t \leq T$,
 \beas
 \hat{Y}_s-Y_s&=&\cE[\hat{Y}_t|\cF_s]-\cE[\xi+zB_\t|\cF_s]=\cE\big[\hat{Y}_t-Y_t
 +\cE[\xi+zB_\t|\cF_t]\big|\cF_s\big]-\cE[\xi+zB_\t|\cF_s]\\
 &=&\cE[\hat{Y}_t-Y_t+\xi+zB_\t|\cF_s]-\cE[\xi+zB_\t|\cF_s]\leq
 \cE^{g_\a}_\g[\hat{Y}_t-Y_t|\cF_s],\q~~ P^\g\mbox{-a.s.}
 \eeas
Applying Theorem \ref{mayaothm} again, we can find an increasing
\cad process $A$ null at $0$ and a process $\bar{Z} \in
\cH^2_\bF([0,T];\hR^d)$ such that
 \beas
 \hat{Y}_t-Y_t=\eta-\xi+\int_t^T \a|\bar{Z}_s|^2ds - A_T+A_t-\int_t^T
 \bar{Z}_s dB^\g_s, \q~~ \fa t \in [0,T],\q P^\g\mbox{-a.s.},
 \eeas
which, in light of the Girsanov Theorem, is equivalent to
 \beas
 \hat{Y}_t-Y_t=\eta-\xi+\int_t^T \neg \big(\a|\bar{Z}_s|^2+\lan \g_s,
 \bar{Z}_s \ran\big)ds - A_T+A_t-\int_t^T \bar{Z}_s dB_s, \q \fa t
 \in [0,T],~~ \pas
 \eeas
On the other hand, we also have
 \beas
 \hat{Y}_t-Y_t=\eta-\xi+\int_t^T (\hat{h}_s-h_s)ds-\int_t^T
 (\hat{Z}_s-Z_s) dB_s, \q~~ \fa t \in [0,T], \q \pas
 \eeas
Thus by comparing the martingale parts and the bounded variation
parts, one has:
 \beas \hat{Z}_t-Z_t \equiv \bar{Z}_t \qq \mbox{and}
 \qq (\hat{h}_t-h_t)dt \equiv \big(\a|\bar{Z}_t|^2+\lan \g_t,
 \bar{Z}_t \ran\big)dt-dA_t,
 \eeas
which implies that $A$ is absolutely continuous and $dA_t=a_tdt$
with $a_t \geq 0$. Consequently,
 \beas
 \hat{h}_t-h_t = \a|\hat{Z}_t-Z_t|^2+\lan \g_t, \hat{Z}_t-Z_t \ran -a_t \le \a|\hat{Z}_t-Z_t|^2+\lan \g_t, \hat{Z}_t-Z_t \ran,
 \q~~ \dtp
 \eeas
This proves the proposition. \qed

We remark that one of the consequences of Proposition \ref{cfm},
especially the representation (\ref{cfm1}), is that the ``\cad
modification" that we found in Proposition \ref{cem} is actually
continuous. In other words, the unique solution of BSDE
\(\ref{BSDEcE}\) should belong to $\hC^\infty_\bF([0,T])$.

We now turn our attention to a comparison theorem for the solutions
to the BSDE (\ref{BSDEcE}). To begin with,
let us note that if $f$ satisfies (H3), then for any $\f \in
L^\infty_\bF([0,T])$, the function $f^\f(t,\o,y) \dfnn
f(t,\o,y)+\f(t,\o)$, $\fa (t,\o,y) \in [0,T] \times \O \times
\hR$, also satisfies (H3). Thus for any $\xi' \in L^\infty(\cF_T)$
and $z \in \hR^d$, the BSDE
 \bea
 \label{fbsde}
 Y_t+zB_t+\int_0^t \big[f(s,Y_s)+\f_s\big]ds= \cE\Big\{\xi' + zB_T+ \int_0^T [
 f(s,Y_s)+\f_s ] ds\Big|\cF_t\Big\}, ~~~ t \in [0,T],
 \eea
admits a unique solution in $\hC^\infty_\bF([0,T])$. We shall
denote this solution by $Y'$.


\begin{thm}\label{compa}
\(Comparison Theorem\) Assume that $f$ satisfies (H3). For fixed $z
\in \hR^d$, let $Y$, $Y' \in \hC^\infty_\bF([0,T])$ be the unique
solution of \(\ref{BSDEcE}\) and \(\ref{fbsde}\) respectively.
Suppose that
 \beas \xi' \geq \xi, \q~~\pas \q\mbox{and}\q \f \geq
0, \q~~ \dtp, \eeas then it holds $P$-a.s. that $Y'_t \geq Y_t$,
$\fa t \in [0,T]$.
\end{thm}

{\it Proof.} We first assume $ \f_t \equiv 0$. For any $\d \in
\hQ^+$, define two stopping times
 \beas
 \si^\d \dfnn \inf\{t \in [0,T)|\;Y'_t \leq Y_t-\d\} \q \mbox{and}
 \q \t^\d \dfnn \inf\{t \in [\si^\d,T]|\;Y'_t \geq Y_t\}.
 \eeas
Here we use the convention that $\inf \es \dfnn T$. Since $Y'_T=\xi'
\geq \xi=Y_T$, $P$-a.s., we must have $\si^\d \leq \t^\d \leq T$,
$P$-a.s.
%
Further, since both $Y$ and $Y'$ have continuous paths, we know that
on $G^\d \dfnn\{\si^\d < T\}$, it holds that
 \bea
 \label{i&ii}
 Y'_{\si^\d} =
 Y_{\si^\d}-\d, \qq Y'_{\t^\d} = Y_{\t^\d}, ~~\pas
 \eea
Next, for a given $t \in [0,T]$, we define a stopping time $\hat{t}
\dfnn t \vee \si^\d \land \t^\d$.
Then, applying Theorem \ref{opts} and
Proposition \ref{psfe}-(iv) we have, $P$-a.s.
 \beas  Y_{\hat{t}}+zB_{\hat{t}}+\int_0^{\hat{t}}
 f(s,Y_s)ds =
 \cE\Big\{ Y_{\t^\d}+zB_{\t^\d}+\int_{\hat{t}}^{\t^\d} \neg
 f(s,Y_s)ds\Big|\cF_{\hat{t}}\Big\}+\int_0^{\hat{t}} f(s,Y_s)ds, \qq \pas
 \eeas
Moreover, since $G^\d \in \cF_{\si^\d}\subset \cF_{\hat{t}}$, we can
deduce from Proposition \ref{psfe} (iii) that
 \bea
 \label{zog} &
 &\b1_{G^\d}\cE\Big\{ \b1_{G^\d} Y_{\t^\d}+zB_{\t^\d}+\neg
 \int_{\hat{t}}^{\t^\d}\neg \neg \b1_{G^\d} f(s,\b1_{G^\d}
 Y_{\hat{s}})ds\Big|\cF_{\hat{t}}\Big\}\nonumber\\
 & =&\b1_{G^\d}\cE\Big\{ \b1_{G^\d} Y_{\t^\d}+zB_{\t^\d}+\neg \int_{\hat{t}}^{\t^\d}\neg
 \neg \b1_{G^\d} f(s, Y_s)ds\Big|\cF_{\hat{t}}\Big\} \\
 &\neg \neg \neg=\neg \neg \neg &\b1_{G^\d}\cE\Big\{
 Y_{\t^\d}+zB_{\t^\d}+\neg \int_{\hat{t}}^{\t^\d}\neg \neg
 f(s,Y_s)ds\Big|\cF_{\hat{t}}\Big\}
 =\b1_{G^\d}Y_{\hat{t}}+\b1_{G^\d}zB_{\hat{t}}. \nonumber
 \eea
By using the $L^\infty$-domination (\ref{infdom}) for $\cE$ and
Proposition \ref{cem} one shows that $P$-a.s.
 \beas
 &&\Big|\cE\Big\{\b1_{G^\d} Y'_{\t^\d}+zB_{\t^\d}+ \neg\int_{\hat{t}}^{\t^\d}\neg
 \neg \b1_{G^\d} f(s,\b1_{G^\d} Y'_{\hat{s}})ds\Big|\cF_r\Big\}
 -\cE\Big\{ \b1_{G^\d} Y_{\t^\d}+zB_{\t^\d}+\neg \int_{\hat{t}}^{\t^\d} \neg \neg
 \b1_{G^\d} f(s,\b1_{G^\d} Y_{\hat{s}})ds\Big|\cF_r\Big\}\Big|\\
 &\neg \neg \neg \le \neg \neg \neg& \Big\| \int_{\hat{t}}^{\t^\d} \neg
 \neg \b1_{G^\d} \big[ f(s,\b1_{G^\d} Y'_{\hat{s}})- f(s,\b1_{G^\d}
 Y_{\hat{s}}) \big] ds \Big\|_\infty \le \k \int_t^T
 \|\b1_{G^\d}Y'_{\hat{s}}- \b1_{G^\d}Y_{\hat{s}}\|_\infty ds,\qq \fa r \in [0,T].
 \eeas
Setting $r=\hat{t}$ in the above and using (\ref{zog}) we obtain
that
 \beas \| \b1_{G^\d} Y'_{\hat{t}}-\b1_{G^\d} Y_{\hat{t}}\|_\infty
 \le \k \int_t^T \|\b1_{G^\d}Y'_{\hat{s}}-
 \b1_{G^\d}Y_{\hat{s}}\|_\infty ds.
 \eeas
The Gronwall inequality then leads to that $\| \b1_{G^\d}
Y'_{\hat{t}}-\b1_{G^\d} Y_{\hat{t}} \|_\infty=0$ for any $t \in
[0,T]$.
In particular, for $t=0$, we obtain that $\b1_{G^\d}
Y'_{\si^\d}=\b1_{G^\d} Y_{\si^\d}$, $P$-a.s., which,  together
with ({\rm i}), shows that $G^\d=\{\si^\d < T\}$ is a null set.
Since $Y'_T \geq Y_T,~P$-a.s. and $\{Y'_t \geq Y_t,\, \fa t \in
[0,T)\}^c \subset \neg \neg \neg \underset{\d \in
\hQ^+}{\bigcup}\neg \neg \{ \si^\d < T\}$,
we conclude that \bea \label{y_y} Y'_t \geq Y_t, \q~~ \fa t \in
[0,T], \q \pas \eea

We now consider the case when $\f_t \geq 0$, $dt \times dP$-a.s.
We proceed as follows. For any $n\in\hN$, let $t^n_j \dfnn
\frac{j}{n}T$, $j=0, 1, \cdot \cdot \cdot, n$ be a partition of
$[0,T]$, and define recursively a sequence of BSDEs:
\beas Y^{j,n}_t+zB_t+\int_0^t f(s,Y^{j,n}_s)ds=
 \cE\Big\{X^n_j+\int_{t_{j-1}^n}^{t_j^n}\f_s
 ds+zB_{t_j^n}+\int_0^{t_j^n}\neg f(s,Y^{j,n}_s)ds\Big|\cF_t\Big\},
 ~~  t \in [0,t_j^n],
 \eeas
where $\{X^n_j\}_{j\ge 0}$ are defined recursively by $X^n_n
=\xi'$, and $X^n_{j-1} \dfnn Y^{j,n}_{t_{j-1}^n}$, for $j=n,\cds,
1$. Now, applying the result for $\f=0$ (similar to (\ref{y_y}))
with $\xi^n_j\dfnn X^n_j+\int_{t^n_{j-1}}^{t^n_j}\f_sds$, we can
then show by induction that for each $1\le j\le n$,
it holds that $Y^{j,n}_t \geq Y_t$, $t \in [0,t_j^n]$, $P$-a.s.
We now define a new process by $Y^n_t \dfnn Y^{j,n}_t$, $t \in
[t_{j-1}^n, t_j^n]$, $j=1, \cdot \cdot \cdot, n$.
It is easy to check that for any $j=1, \cdot \cdot \cdot, n$ and any
$t \in [t^n_{j-1},t^n_j)$,
 \beas Y^n_t+zB_t=\cE\Big\{\xi'+\int_{t_{j-1}^n}^T \f_s ds+zB_T+\int_t^T
 f(s,Y^n_s)ds\Big|\cF_t\Big\},  \qq  \pas
 \eeas
Applying $L^\infty$-domination (\ref{infdom}) for $\cE$ we see
that for any $j=1, \cdot \cdot \cdot, n$ and any $t \in
[t^n_{j-1},t^n_j)$
 \beas
 && \| Y^n_t-Y'_t \|_\infty \\
 & \neg \neg =\neg \neg & \Big\| \cE\Big\{ \xi'+ \neg \int_{t_{j-1}^n}^T \neg \neg
 \f_s ds+zB_T+\neg \int_t^T \neg \neg f(s,Y^n_s)ds\Big|\cF_t\Big\}
 -\cE\Big\{ \xi'+zB_T+ \neg \int_t^T \neg \neg \big[f(s,Y'_s)+\f_s \big]
 ds\Big|\cF_t\Big\}
  \Big\|_\infty \\
 & \neg \neg \leq \neg \neg & \Big\| \int_{t_{j-1}^n}^t \neg \neg
 \f_s ds+ \neg \int_t^T \neg \neg \big(f(s,Y^n_s)- f(s,Y'_s)\big)ds
 \Big\|_\infty \le \frac{T}{n} \| \f \|_\infty + \k \int_t^T
 \neg \neg \|Y^n_s-Y'_s\|_\infty ds.
 \eeas
First applying Gronwall's inequality and then letting $n\to\infty$
we see that
$Y^n_t$ converges to $Y'_t$ in $L^\infty(\cF_t)$, for each $t \in
[0,T]$.
Since both $Y$ and $Y'$ are continuous, we conclude that $Y'_t
\geq Y_t$, $\fa t \in [0,T]$, $P$-a.s. The proof is now complete.
\qed

We can now follow the scheme of \cite{CHMP} and \cite{Pln} to
derive the Doob-Meyer decomposition. For any $Y \in
\hD^\infty_\bF([0,T])$ and $z \in \hR^d$, we define
 \beas
 f^n(t,\o,y) \dfnn n(Y(t,\o)-y), \q~~~ \fa (t,\o,y) \in [0,T]
 \times \O \times \hR, \q~ \fa n \in \hN.
 \eeas
It is easy to check that each $f^n$ satisfies (H3), thus the BSDE
 \bea
 \label{ynf}
 y^n_t+zB_t+\int_0^t f^n(s,y^n_s)ds = \cE\Big\{Y_T + zB_T+\int_0^T
 f^n(s,y^n_s)ds\Big|\cF_t\Big\}, \qq \fa t \in [0,T],
 \eea
admits a unique solution $y^n \in \hC^\infty_\bF([0,T])$. We have
the following lemma.
 \begin{lem}
 \label{complem}
 Assume (H3), and let $y^n$ be the solution of
(\ref{ynf}), $n\ge 1$. Suppose that for a given $Y \in
\hD^\infty_\bF([0,T])$ and $z \in \hR^d$, the process $Y_t+zB_t$, $t
\in [0,T]$ is a $\cE$-submartingale (resp. $\cE$-supermartingale),
then it holds that
 \beas
 y^n_t \geq \(resp.\, \leq\) y^{n+1}_t \geq
 (\mbox{resp. } \leq) Y_t, \q~~ t \in [0,T], \q n \in \hN, \q \pas
  \eeas
\end{lem}

{\it Proof.} We shall prove only the submartingale case, the
supermartingale case is similar. For any $n \in \hN$ and any $\d \in
\hQ^+$, let us define two stopping times
 $$\si^{n,\d} \dfnn \inf\{t
 \in [0,T)|\,y^n_t \leq Y_t-\d\} \q\mbox{ and }\q \t^{n,\d} \dfnn
 \inf\{t \in [\si^{n,\d},T]|\,y^n_t \geq Y_t\}.
 $$
It is easy to see that $\si^{n,\d} \leq \t^{n,\d} \leq T$, $P$-a.s.
Then the right-continuity of $y^n$ and $Y$ leads to that
 \bea
 \label{ynsid}
  y^n_{\si^{n,\d}} \leq Y_{\si^{n,\d}}-\d,~ \; \pas ~\; \mbox{on } \{\si^{n,\d}<T\},
 \q \mbox{and} \q~y^n_{\t^{n,\d}} \ge Y_{\t^{n,\d}}, ~ \; \pas
 \eea
Applying Proposition \ref{psfe}-(iv) and Theorem \ref{opts}, one has
 \beas
 y^n_{\si^{n,\d}}+zB_{\si^{n,\d}}
 =\cE[y^n_{\t^{n,\d}}+zB_{\t^{n,\d}} + \int_{\si^{n,\d}}^{\t^{n,\d}}
 n(Y_s-y^n_s)ds|\cF_{\si^{n,\d}}], \qq \pas
 \eeas
Using (\ref{ynsid}) we deduce that $\int_{\si^{n,\d}}^{\t^{n,\d}}
n(Y_s-y^n_s)ds \geq 0$, $P$-a.s., and combined with
Proposition \ref{psfe}-(i) and Theorem \ref{opts}, we obtain that
 \beas
 y^n_{\si^{n,\d}}+zB_{\si^{n,\d}}
 \geq \cE[Y_{\t^{n,\d}}+zB_{\t^{n,\d}}|\cF_{\si^{n,\d}}] \geq
 Y_{\si^{n,\d}}+zB_{\si^{n,\d}}.
 \eeas
This implies that
 $\{y^n_{\si^{n,\d}} \leq Y_{\si^{n,\d}}-\d\}$ is a null set, thus so is
 $\{\si^{n,\d}<T\}$. Furthermore, since
 \beas \{y^n_t \geq Y_t, ~ t \in [0,T),~ n \in \hN\}^c \subset \neg \neg
 \underset{n \in \hN}{\bigcup}\,\underset{\d \in \hQ^+}{\bigcup}\neg\{
 \si^{n,\d} <T\} \q \mbox{and}\q y^n_T \geq Y_T, ~ n \in \hN,
 \eeas
 it holds $P\{y^n_t \geq Y_t, ~t \in [0,T], ~ n \in \hN\}=1$.
Consequently, we have that $P$-a.s.
 \beas f^n(t, y^n_t)=
 n(Y_t-y^n_t) \geq (n+1)(Y_t-y^n_t)=f^{n+1}(t, y^n_t) ,  \q~~ \fa t
 \in [0,T],\q \fa n \in \hN.
 \eeas
It then follows from Theorem \ref{compa} that $P$-a.s.
 $y^n_t \geq y^{n+1}_t \geq Y_t$, for all $t \in [0,T]$ and $n \in \hN$.
This completes the proof.
\qed

We should note that Lemma \ref{complem} indicates that if $Y_\cd+zB_\cd$ is an
$\cE$-submartingale, then all the processes $\dis A^n_t=\int_0^t
n(y^n_s-Y_s)ds$, $t\ge 0$ are increasing (or decreasing if $Y$ is a
$\cE$-supermartingale),
$\| y^n \|_\infty \leq \|Y \|_\infty \vee \| y^1 \|_\infty$, and
$y^n_t-A^n_t+zB_t$, $t\ge0$ is an $\cE$-martingale. Thus,
Proposition \ref{cfm} implies that there is a unique pair $(h^n,Z^n)
\in L^1_\bF([0,T]) \times \cH^2_\bF([0,T];\hR^d)$ such that
 \bea
 \label{a02}
 y^n_t-A^n_t+zB_t&\neg
 \neg \neg=\neg \neg \neg&y^n_T-A^n_T+zB_T+\int_t^T h^n_s ds
 -\int_t^T Z^n_sdB_s, \q~~t \in [0,T],
 \eea
 and the following estimates hold:
 \bea
 \label{a01}
 -\ell\big(\big|Z^n_t\big|+\big|Z^n_t\big|^2\big)&\neg
 \neg \neg\leq \neg \neg \neg & g_1(t,Z^n_t) \leq h^n_t \leq
 g_2(t,Z^n_t) \leq
 \ell\big(\big|Z^n_t\big|+\big|Z^n_t\big|^2\big), \q  \dtp
 \eea
We shall prove that both $ \{ Z^n \}_{n \in \hN}$ and $\{ A^n_T
\}_{n \in \hN}$ are bounded in a very strong sense.
\begin{lem}\label{ub}
Let the process $Y_t+zB_t$, $t\in [0,T]$, be either an $\cE$-submartingale or an
$\cE$-supermartingale as those in Lemma \ref{complem}, and let
$\{A^n\}$ and $\{Z^n\}$ are processes defined in (\ref{a02}). Then,
for any $p > 0$, $\{ Z^n \}_{n \in \hN}$ is bounded in
$\cM^p(\hR^d)$ and $\{ A^n_T \}_{n \in \hN}$ is bounded in
$L^p(\cF_T)$.
\end{lem}

{\it Proof.} We shall only prove the submartingale case. That is, we
assume that $A^n$ is increasing. From BSDE (\ref{a02}) we see that
 \beas A^n_T &=& y^n_T-y^n_0+\int_0^T h^n_s ds-\int_0^T (Z^n_s-z)
 dB_s,\qq \pas
 \eeas
Let $M \dfnn \| Y \|_\infty \vee \| y^1 \|_\infty$ and use the
domination (\ref{a01}) of $h^n$, we have
 \bea
 \label{estofA}
 |A^n_T|&\leq& 2M+\ell T+ 2\ell \int_0^T
 |Z^n_s|^2ds+\underset{0\leq t \leq T}{\sup}\Big|\int_0^t (Z^n_s-z)
 dB_s\Big|,\q~~ \pas
 \eea
In what follows for each $p>0$ we denote $C_p>0$ to be a generic
constant depending only on $p$, as well as $\ell,T,M,|z|$, which may
vary from line to line.
Using (\ref{estofA}) and the Burkholder-Davis-Gundy inequality one
shows that
 \beas
 E|A^n_T|^p \le C_p\Big\{1\neg +\neg E\Big[\int_0^T
 \neg|Z^n_s|^2ds\Big]^p\neg +\neg E\Big[\int_0^T \neg
 |Z^n_s-z|^2ds\Big]^{p/2}\Big\} \le C_p\Big\{1\neg +\neg E\Big[\neg\int_0^T \neg
 |Z^n_s-z|^2ds\Big]^p\Big\}.
 \eeas
Thus it suffices to show that $\underset{n \in \hN}{\sup}
E\Big(\int_0^T |Z^n_s-z|^2ds\Big)^p < \infty$. For any $\a>0$, we
apply It\^o's formula to $\dis e^{\a y^n_t}$ to get:
 \bea
 \label{Ito}
 & & e^{\a y^n_0}+\frac{\a^2}{2}\int_0^T e^{\a y^n_s}|Z^n_s-z|^2 ds \nonumber\\
 &=& e^{\a y^n_T}+\a  \Big[\int_0^T e^{\a y^n_s} h^n_s ds
 -\int_0^T e^{\a y^n_s}d A^n_s-\int_0^T e^{\a y^n_s}(Z^n_s-z)dB_s\Big] \\
 &\leq & e^{\a y^n_T}+\a \ell  \int_0^T e^{\a y^n_s}ds+ 4\a\ell \int_0^T e^{\a y^n_s}
 |Z^n_s-z|^2ds+4\a \ell \int_0^T e^{\a y^n_s}|z|^2ds\nonumber\\
 & & -\a \int_0^T e^{\a y^n_s}(Z^n_s-z)dB_s.\nonumber
 \eea
Note that the last inequality is due to the fact that $A^n$ is
increasing. It then follows that
 \beas
 (\frac{\a^2}{2}-4\a\ell)\int_0^T e^{\a y^n_s}|Z^n_s-z|^2 ds \le
 C_p+\a\, \underset{0\leq t \leq T}{\sup}\Big|\int_0^t
 e^{\a y^n_s}(Z^n_s-z)dB_s\Big|.
 \eeas
Choose $\a>8\ell$, and applying the Burkholder-Davis-Gundy
inequality again we obtain that
 \beas
 && E\Big(\int_0^T \neg \neg e^{\a y^n_s}|Z^n_s-z|^2 ds \Big)^p \leq C_p+ C_p
 E\Big(\int_0^T \neg \neg e^{2 \a y^n_s}|Z^n_s-z|^2ds\Big)^{p/2}\\
 &\leq& C_p+C_p e^{pM\a/2}E\Big(\int_0^T \neg \neg e^{\a
 y^n_s}|Z^n_s-z|^2ds\Big)^{p/2}\leq C_p+\frac{1}{2}E\Big(\int_0^T
 \neg \neg e^{\a y^n_s}|Z^n_s-z|^2ds\Big)^p,
 \eeas
which implies that $E\Big(\int_0^T \neg \neg e^{\a y^n_s}|Z^n_s-z|^2
ds \Big)^p$ is dominated by a constant independent of $n$. This
proves the lemma in the submartingale case. The supermartingale case
can be proved in the same way except that in (\ref{Ito}) the It\^o's
formula should be applied to  $e^{-\a y^n_t}$. The proof is now
complete. \qed

We are now ready to prove the Doob-Meyer Decomposition Theorem.

\begin{thm}
 \label{dqfm}
Assume that $\cE$ is a regular quadratic $\cF$-expectation
satisfying the one-sided $g$-domination (\ref{1sdom}). For any $Y
\in \hC^\infty_\bF([0,T])$ and any $z \in \hR^d$, if the process
$Y_t+zB_t$, $t \in [0,T]$, is an $\cE$-submartingale (resp.
$\cE$-supermartingale), then there exists a continuous increasing
(resp. decreasing) process $A$ null at $0$ such that $Y_t-A_t+zB_t$,
$t\ge 0$, is a local $\cE$-martingale. Furthermore, if $A$ is
bounded, then $Y_t-A_t+zB_t$, $t\ge 0$, is an $\cE$-martingale.
\end{thm}

{\it Proof.} We again prove only the submartingale case, as the
submartingale case is similar. To begin with, let $y^n$ be the
solutions to (\ref{ynf}), $n=1,2,\cds$, and still denote $M \dfnn \|
Y \|_\infty \vee \| y^1 \|_\infty$. Since $y^n\ge Y$, by the
definition of processes $A_n$'s and Lemma \ref{complem}, we see that
 $$ E\int_0^T |y^n_s-Y_s|ds = \frac{1}{n}E[|A^n_T|] \leq
 \frac{1}{n}\underset{n \in \hN}{\sup}\| A^n_T \|_1 \to 0,
 $$
as $n\to\infty$. Moreover, since $y^n$'s converges decreasingly to
$Y$, and $Y$ is continuous, we can further conclude, in light of
Dini's Theorem,  that $\pas$
 \bea
 \label{dini}
 \lim_{n\to\infty}\underset{t \in [0,T]}{\sup}(y^n_t-Y_t) =0,
\q  \mbox{thus} \q \lim_{m,n\to\infty}\underset{t \in
[0,T]}{\sup}|y^m_t-y^n_t| =0.
 \eea

We first show that there exists a subsequence of $\{A^n\}$, still
denoted by $\{A^n\}$, such that the sequence $\{A^n_T\}_{n \in \hN}$
is uniformly integrable. To see this, we claim that the processes $Z^n$ converges to some
process  $Z$ in $\cH^2_\bF([0,T];\hR^d)$, as $n\to\infty$. In fact,
applying It\^o's formula to $|y^m_t-y^n_t|^2$ on $[0,T]$ we obtain
\bea
 &&|y^m_0-y^n_0|^2+\int_0^T |Z^m_s-Z^n_s|^2 ds \label{zconv} \\
 &\tneg \dneg =\dneg \tneg &|y^m_T- y^n_T|^2\neg+\neg 2\neg
 \int_0^T \dneg (y^m_s\neg-\neg y^n_s)\big[(h^m_s-h^n_s)ds\neg-\neg
 (dA^m_s\neg-\neg dA^n_s)\neg-\neg
 (Z_s^m\neg-Z^n_s)dB_s \big] \nonumber\\
 &\tneg \dneg \le \dneg \tneg & |y^m_T- y^n_T|^2\dneg+\neg 2
 \dneg\underset{s \in
 [0,T]}{\sup}\neg|y^m_s-y^n_s|\Big\{\dneg\int_0^T
 \dneg 2\ell\big(1\neg+\neg |Z^m_s|^2+\neg |Z^n_s|^2\big)ds\neg+\neg A^m_T
 +\neg A^n_T\Big\} \nonumber \\
 &&-\neg 2 \int_0^T \dneg
 (y^m_s\neg-\neg y^n_s)(Z_s^m\neg-\neg Z^n_s)dB_s. \nonumber
 \eea
Taking expectation on both sides of (\ref{zconv}) and applying
H\"older's inequality one has
 \beas
& &\tneg \dneg E\Big\{\int_0^T |Z^m_s-Z^n_s|^2 ds\Big\}\\
&\tneg \tneg \le& \tneg \dneg E\Big\{\underset{s \in
[0,T]}{\sup}|y^m_s\neg-\neg y^n_s|^2\neg\Big\}\neg+\neg 2\bigg\{\neg
E\Big[\underset{s \in [0,T]}{\sup}|y^m_s\neg-\neg y^n_s|^2\Big]\neg
E\Big[\neg\int_0^T \dneg 2\ell\big(1\neg+\neg |Z^m_s|^2\neg+\neg
|Z^n_s|^2\big)ds\neg+\neg A^m_T
\neg +\neg A^n_T\Big]^2\neg \bigg\}^{1/2}\\
 &\tneg \tneg \le &\tneg \dneg  E\Big\{\neg\underset{s \in [0,T]}{\sup}|y^m_s\neg-\neg y^n_s|^2\neg\Big\}
 \neg+\neg C\bigg\{\neg E\Big[\underset{s \in [0,T]}{\sup}|y^m_s\neg-\neg y^n_s|^2\Big]\bigg\}^{1/2}
  \Big[1\neg+\neg \underset{k \in \hN}{\sup}\|Z^k\|^2_{\cM^4}\neg+\neg \underset{k \in
  \hN}{\sup}\|A^k_T\|_{L^2(\cF_T)}\Big],
 \eeas
 where $C>0$ is a constant depending only on $\ell$ and $T$. This, together with Lemma \ref{ub},
imply that $\{Z^n\}_{n \in \hN}$ is a Cauchy sequence in
$\cH^2_\bF([0,T];\hR^d)$, hence has a limit $Z\in
\cH^2_\bF([0,T];\hR^d)$. A simple application of the
Burkholder-Davis-Gundy inequality leads to that
 \bea
 \label{f01}
 \underset{t\in [0,T]}{\sup}\Big|\int_0^t (Z^n_s-Z_s)
 dB_s\Big|\to 0 ~~ \mbox{in}~~ L^2(\cF_T),~~ \mbox{as }\, n
 \to \infty.
 \eea
Applying \cite[Lemma 2.5]{Ko} we can find a subsequence of
$\{Z^n\}_{n \in \hN}$, still denoted by $\{Z^n\}_{n \in \hN}$, such
that
$\underset{n}{\sup}|Z^n| \in \cH^2_\bF([0,T];\hR^d)$ and that
$\underset{n}{\sup}|\int_0^T (Z^n_s-z) dB_s| \in L^2(\cF_T)$. Then
in light of (\ref{a01}) and (\ref{a02}), it holds $P$-a.s. that for
any $n \in \hN$
 \beas
 A^n_T &=& y^n_T-y^n_0+\int_0^T h^n_s ds
 -\int_0^T (Z^n_s-z) dB_s \\
 &\le& 2M+\ell T +2\ell \int_0^T \underset{n}{\sup}|Z^n_s|^2 ds+
 \underset{n}{\sup}\Big|\int_0^T (Z^n_s-z) dB_s\Big|  \in L^1(\cF_T).
 \eeas
We can then deduce that $\dis \underset{n \in \hN}{\sup}A^n_T \in
L^1(\cF_T)$, which implies that,  $P$-almost surely, $A^n_t
\le E\big[ \underset{n \in \hN}{\sup}A^n_T\big|\cF_t\big]$, for all $t \in [0,T]$,
$n\in \hN$. Now let us define a sequence of stopping times
 \bea
 \label{tkm}
 \t_{k} \dfnn \inf\{t \in [0,T]: E\big[ \underset{n \in \hN}{\sup}A^n_T\big|\cF_t\big] > k \}\land
 T,\qq k \in \hN.
 \eea
Clearly,
$\t_{k} \nearrow T$, $\pas$, as $k \to \infty$. Furthermore, let us
denote
$p_k \dfnn p (k+M,|z|)$, $J_k \dfnn J (k+M,|z|)$ and $\a_k \dfnn \a
(k+M,|z|)$, and define $Y^k_t\dfnn Y_{t \land \t_k}$, $y^{n,k}_t
\dfnn y^n_{t \land \t_k}$, $A^{n,k}_t \dfnn A^n_{t \land \t_k}$, $
t\in[0,T]$. We will show that for any $k \in \hN$, there exists a
subsequence of $\{A^n\}_{n \in \hN}$, denoted again by $\{A^n\}_{n
\in \hN}$ itself, such that for all $k\in\hN$, it holds that
$\underset{n \to \infty}{\lim}A^{n,k}_{t}=\tilde{A}^k_t$, $t \in
[0,T]$, $\pas$ for some continuous, increasing process
$\tilde{A}^k$.

To see this, let us first fix $k \in \hN$. For each $n\in\hN$,
applying Theorem \ref{opts} and Proposition \ref{cem} we have
 \beas y^{n,k}_t-A^{n,k}_t+zB_{t
 \land \t_k}=\cE[y^n_{\t_k}-A^n_{\t_k}+zB_{\t_k}|\cF_t], \qq \fa t
 \in [0,T].
 \eeas
Applying Proposition \ref{cfm}, we can find a unique pair
$(h^{n,k},Z^{n,k}) \in L^1_\bF([0,T]) \times \cH^2_\bF([0,T];\hR^d)$
such that \bea
 \label{b02}
 && y^{n,k}_t\neg-\neg
 A^{n,k}_t=y^{n,k}_T \neg - \neg A^{n,k}_T+\neg \int_t^T \neg \neg
  h^{n,k}_s ds - \neg\int_t^T \neg \neg \big(Z^{n,k}_s-\b1_{\{s \leq
 \t_k\}}z\big) dB_s, \q \fa t \in [0,T].
 \eea
On the other hand, by (\ref{a02}) we have
 \bea\label{b03}
 y^{n,k}_t\neg-\dneg A^{n,k}_t\neg=\neg y^{n,k}_T \dneg -\dneg A^{n,k}_T\neg + \neg
 \int_t^T \b1_{\{s \le \t_k\}}
 h^n_s ds\neg - \neg\int_t^T \b1_{\{s \le \t_k\}} (Z^n_s-z) dB_s,\q \fa t \in [0,T].
 \eea
Thus by comparing the martingale parts and the bounded variation
parts of (\ref{b02}) and (\ref{b03}), one has $h^{n,k}_t \equiv
\b1_{\{t \le \t_k\}}h^n_t$ and $Z^{n,k}_t\equiv \b1_{\{t \le
\t_k\}}Z^n_t$. Moreover, it also follows from Lemma \ref{cfm} that
there is a BMO process $\g^{n,k}$ with $\| \g^{n,k} \|^2_{BMO}\leq
J_k$ such that
 \bea
 \label{b01}
  && -\a_k|Z^{m,k}_t\neg-\neg
Z^{n,k}_t|^2\neg +\neg \lan \g^{m,k}_t,Z^{m,k}_t\neg-\neg Z^{n,k}_t
\ran \le h^{m,k}_t \neg- \neg h^{n,k}_t \nonumber \\
&&\qq\qq\qq\qq\leq  \a_k|Z^{m,k}_t \neg- \neg Z^{n,k}_t|^2 +\lan
\g^{n,k}_t,Z^{m,k}_t\neg- \neg Z^{n,k}_t \ran, \qq \dtp \eea

Note that (\ref{b01}) implies that for any $m,n \in \hN$,
 \beas
 && E\int_0^{\t_k} \neg \neg \big|h^m_s\neg-\neg h^n_s\big|ds
  \leq E\int_0^{\t_k} \neg \neg \Big[ \a_k\big|Z^m_s\neg-\neg Z^n_s\big|^2 \neg
 +\neg \big(\big|\g^{m,k}_s\big|\vee\big|\g^{n,k}_s\big|\big)\big|Z^m_s\neg-\neg Z^n_s \big|\Big]ds \\
 &  \neg \neg\le  \neg \neg & \a_k E\neg
 \int_0^T \neg \neg \big|Z^m_s\neg -\neg Z^n_s\big|^2 ds
 \neg+\neg \bigg\{E\neg \int_0^T \neg\neg
 \big(\big|\g^{m,k}_s\big|^2 \neg+\neg
 \big|\g^{n,k}_s\big|^2\big)ds\, E\neg \int_0^T \neg\neg
 \big|Z^m_s\neg-\neg Z^n_s \big|^2ds
 \bigg\}^{\frac12}.
 \eeas
Hence, one can deduce from the convergence of $Z^n$ in
$\cH^2_\bF([0,T];\hR^d)$ that $ \big\{\b1_{\{\cdot \land
\t_k\}}h^n_\cdot \big\}_{n \in \hN}$ is a Cauchy sequence in
$L^1_\bF([0,T])$.
Let $\tilde{h}^k$ be its limit in $L^1_\bF([0,T])$, 
it then follows that
\bea\label{f02} \hb{$\dis \underset{t\in
[0,T]}{\sup}\Big|\int_0^{t\land \t_k}(h^n_s-\tilde{h}^k_s)ds\Big|
\to 0~$ in $~L^2(\cF_{\t_k})$,\; as $\,n \to \infty$.} \eea

Now let us define $\tilde{A}^k_t \dfnn Y^k_t- Y^k_0+\int_0^{t \land
\t_k} \tilde{h}^k_s ds -\int_0^{t \land \t_k}(Z_s-z) dB_s$, $ t \in
[0,T]$. Clearly, $\tilde A^k$ is continuous. Furthermore, since
 \beas
 A^{n,k}_t=y^{n,k}_t-y^{n,k}_0+\int_0^{t \land \t_k} h^n_s ds -\int_0^{t \land \t_k}
 (Z^n_s-z) dB_s, \q~~ \fa t \in
 [0,T], \q \fa n \in \hN,
 \eeas
applying Bounded Convergence Theorem as well as (\ref{dini}),
(\ref{f01}) and (\ref{f02}), one shows that
$\dis \underset{t \in [0,T]}{\sup} |A^{n,k}_t-\tilde{A}^k_t|$
converges to $0$ in $L^1(\cF_{\t_k})$, 
as $ n \to\infty$. Therefore, we can find a subsequence of
$\{A^n\}_{n \in \hN}$, still denoted by $\{A^n\}_{n \in \hN}$, such
that
 \bea
 \label{f03}
 \underset{n \to \infty}{\lim}A^{n,k}_t=\tilde{A}^k_t, \q~~~ \fa t
 \in [0,T], \q~\pas \eea
We note that (\ref{f03}) indicates that $\tilde{A}^k$ is an
increasing process. Furthermore, applying the Helly Selection
Theorem if necessary, we can assume that the convergence in
(\ref{f03}) holds true for all $k\in\hN$ for this subsequence.

We can now complete the proof. By the definition of $\t_k$
(\ref{tkm}) and the continuity of $A^n$, one can deduce that for any
$k$, $n \in \hN$, $A^n_{\t_k}\le k$, $\pas$

Hence for any $k \in \hN$, (\ref{f03}) implies that $\pas$
 \beas
 |\tilde{A}^k_t| \le k, \q~ \fa t \in [0,T] \q~ \mbox{and}\q~
 \tilde{A}^k_t \equiv \tilde{A}^k_{\t_k},\q~ \fa t \in [\t_k,T].
 \eeas
Note that $\tilde{A}^k_t=\underset{n \to
\infty}{\lim}A^{n,k}_t=\underset{n \to
\infty}{\lim}A^{n,k+1}_{t\land \t_k}=A^{k+1}_{t\land \t_k}$, $t \in
[0,T]$, $\pas$ we can define a continuous, increasing process $A_t
\dfnn \tilde{A}^k_t$, $t \in [0,\t_k]$, $k \in \hN$. Clearly, $A$ is
null at $0$. For fixed $k \in \hN$, and $t \in [0,T]$, applying the
$L^{p_k}$-domination (\ref{pdom}) of $\cE$ yields that
 \beas \|
 \cE[y^n_{\t_k}-A^n_{\t_k}+zB_{\t_k}|\cF_t]-\cE[Y_{\t_k}-A_{\t_k}+zB_{\t_k}|\cF_t]
 \|_{p_k} \leq 3\| y^n_{\t_k}-Y_{\t_k}\|_{p_k}+3\|
 A^n_{\t_k}-A_{\t_k}\|_{p_k}.
 \eeas
By considering a subsequence, we have, $P$-a.s.
 \beas Y_{t \land
 \t_k}- A_{t \land \t_k} + zB_{t \land \t_k}
 &=&\underset{n \to \infty}{\lim}\big(y^{n,k}_{t \land \t_k} - A^{n,k}_{t \land \t_k}+
 zB_{t \land \t_k}\big)\\
 &=& \underset{n \to \infty}{\lim} \cE[y^{n,k}_{\t_k}- A^{n,k}_{\t_k} +
 zB_{\t_k}|\cF_t]= \cE[Y_{\t_k} - A_{\t_k}+ zB_{\t_k}|\cF_t].
 \eeas
Then, Proposition \ref{cem}, together with the continuity of $Y$ and
$A$, implies that $P$-a.s.
 \bea
 \label{locm} Y_{t \land \t_k}-A_{t
 \land \t_k}+zB_{t \land
 \t_k}=\cE[Y_{\t_k}-A_{\t_k}+zB_{\t_k}|\cF_t],\qq \fa t \in [0,T].
 \eea
In other words, $Y_t-A_t+zB_t$, $t\ge 0$ is a local
$\cE$-martingale, proving the first part of the theorem.

To see the last part of the theorem, we assume further that $A$ is
bounded. Let $K=\| Y \|_\infty+\| A \|_\infty$, $R=|z|$ and
$p=p(K,R)$. Fix a $t \in [0,T]$, applying $L^p$-domination
(\ref{pdom}) again we obtain that for any $k \in \hN$,
 \beas
 &&\| \cE[Y_{\t_k}- A_{\t_k}+
 zB_{\t_k}|\cF_t]-\cE[Y_T- A_T+zB_T|\cF_t]
 \|_p\\
 &&\qq\qq\le  R \| B_{t \land \t_k}- B_t \|_p+ 3\|Y_{\t_k}-Y_T \|_p +
  3\|A_{\t_k}- A_T \|_p+C_R\|T-\t_k\|_p.
 \eeas
Clearly, the right hand side above converges to $0$ as $k \to
\infty$, thanks to the Burkholder-Davis-Gundy inequality and the
Bounded Convergence Theorem. Thus, taking a subsequence if
necessary, we may assume that
$\cE[Y_{\t_k}-A_{\t_k}+zB_{\t_k}|\cF_t]$ converges $P$-a.s. to
$\cE[Y_T-A_T+zB_T|\cF_t]$.
Letting $k \to \infty$ in (\ref{locm}), the continuity of $Y$ and
$A$ imply that \beas Y_t-A_t+zB_t=\cE[Y_T-A_T+zB_T|\cF_t], \qq \pas
\eeas Eventually, applying Proposition \ref{cem} and using the
continuity of $Y$ and $A$ again we have $P$-a.s.
 \beas
 Y_t-A_t+zB_t=\cE[Y_T-A_T+zB_T|\cF_t], \qq \fa t \in [0,T],
 \eeas
which means that $Y_t-A_t+zB_t$, $t\ge0$, is an $\cE$-martingale. The
proof is now complete. \qed

\section{Representation Theorem of Quadratic $\cF$-Expectations}
\setcounter{equation}{0}

In this section we prove the representation theorem for quadratic
$\cF$-expectations. We assume that $\cE$ is a translation
invariant quadratic $\cF$-expectation dominated by two quadratic
$g$-expectations $\cE^{g_1}$ and $\cE^{g_2}$ from below and above,
and both $g_1$ and $g_2$ satisfy (H1) and (H2) with the same
constant $\ell>0$. We also assume that $\cE$ satisfies the
$L^p$-domination (\ref{pdom}), the $L^\infty$-domination
(\ref{infdom}), and the one-sided $g$-domination (\ref{1sdom}).

We begin our discussion by considering the following special
semi-martingale:
 \bea
 \label{Yz}
 Y^z_t \dfnn \ell(|z|+|z|^2)t+zB_t,\qq \fa t \in [0,T], \qq z \in
 \hR^d.
 \eea
By the comparison theorem of BSDEs, it is easy to see that $Y^z$
is an $\cE^{g_1}$-submartingale, whence an $\cE$-submartingale.
Then, by the Doob-Meyer decomposition (Theorem \ref{dqfm}) there
exists a continuous, increasing process $A^z$ null at $0$ such
that $Y^z-A^z$ is a local $\cE$-martingale. We claim that
$A^z_T\in L^\infty(\O)$, and hence $Y^z-A^z$ is a true
$\cE$-martingale. Indeed, let $\{\t^z_n\}_{n\ge1}$ be a sequence
of ``reducing" stopping times, that is, $\t^z_n \nearrow T$,
$\pas$, such that
 \bea
 \label{dqfm1}
 Y^{z,n}_t-A^{z,n}_t=\cE[Y^{z,n}_T-A^{z,n}_T|\cF_t],\q~~ \fa t \in
 [0,T],\q \pas,
 \eea
where $Y^{z,n}_t \dfnn Y^z_{t \land \t^z_n}$, $ A^{z,n}_t \dfnn
A^z_{t \land \t^z_n}$, $\fa t \in [0,T]$. For any $n \in \hN$, we
know from Proposition \ref{cfm} that there is a unique pair
$(h^{z,n},Z^{z,n}) \in L^1_\bF([0,T]) \times \cH^2_\bF([0,T];\hR^d)$
such that
 \bea
 \label{c02}
 &&Y^{z,n}_t=Y^{z,n}_T-A^{z,n}_T+A^{z,n}_t+\int_t^T
 h^{z,n}_sds -\int_t^T Z^{z,n}_sdB_s, \q \fa t \in [0,T],
 \eea
such that the generator $h$ satisfies the following estimate:
 \bea
 \label{c01}
 &&\hspace{-1cm}-\ell\big(\big|Z^{z,n}_t\big|\neg+\neg
 \big|Z^{z,n}_t\big|^2\big)\neg \leq \neg g_1(t,Z^{z,n}_t) \leq
 h^{z,n}_t \neg \leq \neg g_2(t,Z^{z,n}_t)\neg \leq \neg
 \ell\big(\big|Z^{z,n}_t\big|\neg +\neg \big|Z^{z,n}_t\big|^2\big),
 \q  \dtp
 \eea
Comparing (\ref{Yz}) and (\ref{c02}) we see that
 \bea
 \label{AZeq}
 dA^{z,n}_t-h^{z,n}_tdt\equiv \b1_{\{t \leq
 \t^z_n\}}\ell(|z|+|z|^2)dt \q \mbox{and} \q  Z^{z,n}_t \equiv
 \b1_{\{t \leq \t^z_n\}}z.
 \eea
This, together with (\ref{c01}), implies that $P$-a.s.
 \beas
 A^{z,n}_T=\int_0^T h^{z,n}_tdt+ \int_0^T \b1_{\{t \leq
 \t^z_n\}}\ell(|z|+|z|^2)dt \le 2\ell(|z|+|z|^2)T.
 \eeas
Letting $n \to \infty$ we obtain that $A^z_T$ is bounded by
$2\ell(|z|+|z|^2)T$,
proving the claim.

Now, in light of  Proposition \ref{cfm}, we can assume that there
exists a unique pair $(h^z,Z^z) \in L^1_\bF([0,T]) \times
\cH^2_\bF([0,T];\hR^d)$ such that (\ref{c02})---(\ref{AZeq}) hold.
In other words,
denoting
 \bea\label{repg} g(t,\o,z)\dfnn h^z_t(\o), \qq
 (t,\o,z)\in[0,T]\times \O \times \hR^d,
 \eea
it holds that
 \bea
 \label{e02}
 Y^z_t-A^z_t& = & Y^z_T-A^z_T+\int_t^T
 g(s,z)ds-\int_t^T zdB_s, \qq  t \in [0,T],\\
 \label{e01}
  -\ell\big(|z|+|z|^2\big)&\neg \neg \le \neg \neg& g_1(t,z) \le
 g(t,z) \le g_2(t,z) \le \ell\big(|z|+|z|^2\big), \q\dtp, \bs\\
 \label{e03}
 dA^z_t&=&g(t,z)dt +\ell(|z|+|z|^2)dt, \qq t\in[0,T].
 \eea
We shall show that $g$ is the desired representation generator of
the quadratic $\cF$-expectation $\cE$.

To begin with, let us define, for any $z,z'\in\hR^d$, a function
 \bea
 \label{gmu}
 g^{z,z'}_\ell(v)\dfnn \ell(1+|z|+|z'|)|v|, \qq \forall v\in\hR^d,
 \eea
and denote the corresponding $g^{z,z'}_\ell$-expectation by
$\cE^{z,z'}_\ell(\cd)$. It is worth noting that
$\cE^{z,z'}_\ell(\cd)$ is a Lipschitz $g$-expectation studied in
\cite{BCHMP} and \cite{Pln}. We should note here that if $g$ is a
quadratic generator
satisfying (H1) and (H2), then it must satisfy a ``local Lipschitz
property" which can be written as
 \bea
 \label{locallip}
 |g(t,z)-g(t,z')|\le \ell(1+|z|+|z'|)|z-z'|=g^{z,z'}_\ell(|z-z'|),
 \qq\forall z,z'\in\hR^d.
 \eea

Now let $g$ be a given deterministic quadratic generator satisfying (H1) and (H2).
For fixed $z\in\hR^d$, consider the process $Y^{g,z}_t\dfnn
\cE^g\{zB_T|\cF_t\}$, $t\ge 0$. Since $zB_T\in \sL^\infty_T$, we
know that (recall the BSDEs (\ref{qBSDE}) and (\ref{qBSDE1}))
$Y^{g,z}_t$ must have the following explicit expression:
 \bea \label{Ygz} Y^{g,z}_t=zB_t+\int_t^T
 g(s,z)ds, \qq t\in[0,T].
 \eea

\if{0} We should note that this $g$-expectation is defined in a
generalized sense, because the terminal value $zB_T$ is not bounded,
but rather one that has |z|^2T}<\infty$.)

\begin{lem}
\label{Ygzform}
Suppose that $g$ is a deterministic function
satisfying (H1) and (H2). Then $Y^{g,z}_t$ has the following
explicit expression: \bea \label{Ygz} Y^{g,z}_t=zB_t+\int_t^T
g(s,z)ds, \qq t\in[0,T]. \eea
\end{lem}

{\it Proof.} Note that the main difficulty of this result comes from
the unboundedness of $B_T$, and the quadratic growth of $g$. To
overcome this difficulty we recall the definition of
$\cE^g\{zB_T|\cF_t\}$. That is, $Y^{g,z}_t$ is the monotone limit of
the sequence $\{Y^{g,z,n}_t\}_{n \ge 1}$, where each $Y^{g,z,n}$ is
the unique solution to the BSDE:
\bea\label{Ygzn} Y^{g,z,n}_t=zB_{T}+\int_{t}^{T}
g_n(s,Z^{g,z,n}_s)ds-\int_{t}^{T}Z^{g,z,n}_sdB_s, \qq t\in[0,T],
\eea
where $g_n$'s are Lipschitz in $z$ and $g_n\ua g$. Note that if $g$
is deterministic, then so are $g_n$'s. But for the deterministic
function $g_n$, the process
$$Y_t\dfnn zB_t+\int_t^T g_n(s,z)ds, \qq t\ge 0
$$
is an adapted solution to the BSDE (\ref{Ygzn}). Thus by uniqueness
$Y^{g,z,n}_t\equiv Y_t=zB_t+\int_t^T g_n(s,z)ds$, $t\ge 0$, for all
$n$. The Lemma then follows by letting $n\to\infty$. \qed \fi

Let us fix $z,z'\in \hR^d$, and define
%
%
$\label{hatE} \hat\cE^{z,z'}_t\dfnn Y^{g,z}_t-Y^{g,z'}_t=
(z-z')B_t+\int_t^T(g(s,z)-g(s,z'))ds$, $t\ge 0$. We have the
following lemma.
\begin{lem}
\label{hatEmu} Assume that $g$ is a deterministic function satisfying
(H1) and (H2). Then the process $\xi_t\dfnn\hat{\cE}^{z,z'}_t$, $t\ge 0$
is a $\cE^{z,z'}_\ell$-submartingale.
\end{lem}

{\it Proof.}
%
For any $s\le t$, define
\bea \label{tildeY} \tilde{Y}_s&\dfnn& \cE^{z,z'}_\ell\{\hat
\cE^{z,z'}_t|\cF_s\}=\Big[(z-z')B_t+\int_t^T(g(r,z)-g(r,z'))dr\Big]\nonumber\\
&&+\int_s^t\m(1+|z|+|z'|)|\tilde Z_r|dr-\int_s^t\tilde Z_rdB_r. \eea
Since $g$ is deterministic, the BSDE (\ref{tildeY}) has a unique
solution $(\hat Y,\hat Z)$, where
$$\hat Y_s\dfnn
(z-z')B_s+\int_t^T(g(r,z)-g(r,z'))dr+\int_s^t\m(1+|z|+|z'|)|z-z'|dr,
$$
and $\hat Z\equiv z-z'$.  Thus, denoting $\d g(r)\dfnn
g(r,z)-g(r,z')$, we have
\beas\tilde Y_s=\hat Y_s &=& (z-z')B_s+\int_t^T\d
g(r)dr+\int^t_s\m(1+|z|+|z'|)|z-z'|dr\\
&=&(z-z')B_s+\int_s^T\d g(r)dr+\int^t_s\{\m(1+|z|+|z'|)|z-z'|-\d
g(r)\}dr \\
&\ge & (z-z')B_s+\int_s^T \d g(r)dr. \eeas
But by definition of $\hat\cE^{z,z'}$ we see that the right hand
side above is exactly $
\hat\cE^{z,z'}_s=\xi_s$. This, combined with (\ref{tildeY}), shows
that $\xi=\hat\cE^{z,z'}$ is an $\cE^{z,z'}_\ell$-submartingale.
\qed


We now introduce  some extra assumptions on the quadratic
$\cF$-expectation $\cE$, which will be useful in the study of the
representation theorem. The first one is motivated by Lemma
\ref{hatEmu}.
\begin{itemize}
\item[{\bf (H4)}] There exists a constant $\m>0$, such that for
any fixed $z, z'$, it holds that
 \bea
 \label{Edomin}
 \cE\{z B_T|\cF_t\}-\cE\{z'B_T|\cF_t\}\le \cE^{z,z'}_\m\{(z-z')B_T|\cF_t\}.
 \eea
\end{itemize}

The next assumption extends the ``translation invariance" of the
nonliear expectation $\cE$.
 \begin{itemize}
  \item[{\bf (H5)}] For any $z \in \hR^d$, $\t \in \cM_{0,T}$, $0\le t < \tilde{t}\le
T$, and $\xi \in L^\infty(\cF_{\tilde{t} \land \t})$, it holds
that
 \bea\label{Transinv}
 \cE[\xi+zB_{\tilde{t} \land \t}-zB_{t \land \t}|\cF_t]=\cE[\xi+zB_{\tilde{t} \land
 \t}|\cF_t]-zB_{t \land \t},\q \pas
 \eea
\end{itemize}


We note that the assumption (H5) is not a consequence of
Proposition \ref{psfe}-(iv), since the random variable $zB_t$ is
not bounded(!).
However, the left hand side of (\ref{Transinv}) is well defined,
since
$\xi+zB_{\tilde{t} \land \t}-zB_{t \land \t}=\xi+\int_t^{\tilde
t}z{\bf 1}_{\{s\le \t\}}dB_s\in\tilde{\sL}^\infty_T\subset$
Dom$(\cE)$ (see Remark \ref{extLinf}).

Finally, we give an assumption that essentially states that the
process $\{zB_t\}_{t\ge0}$ has the ``independent increments"
property under the nonlinear expectation $\cE$.
 \begin{itemize}
 \item[{\bf (H6)}] ~For any $z \in \hR^d$, and any $0\le s\le t\le T$,
 it holds that
 \bea
 \label{indep}
  \cE[z(B_t-B_s)|{\cal F}_s]=\cE[z(B_t-B_s)],\q~~ \pas
 \eea
 \end{itemize}
The following Lemma is more or less motivated the assumption (H6),
and it will play an important role in the proof of the
representation theorem.
 \begin{lem}
 \label{dll}
 Assume that $\cE$ is a regular quadratic $\cF$-expectation satisfying
 (H6).
Then the random function $g$ defined in (\ref{repg}) is
deterministic, and it holds that
 \bea
 \label{limg}  g(t,z)=\lim_{h \to 0} \frac{\cE\{z(B_{t+h}-B_t)\}}{h},\q~~
 \pas, \q \fa (t,z) \in [0,T] \times \hR^d.
 \eea
Moreover, if in addition $\cE$ satisfies (H4), then $g$ is local
Lipschitz continuous.
\end{lem}

{\it Proof.}
We first show that $g$ is deterministic. To this end, we fix $z
\in \hR^d$. For any $0\le t < t+h\le T$, one can deduce from
(\ref{e02}) that
 \beas
 z(B_{t+h}-zB_{t})-\int_t^{t+h}g(s,z)ds=Y^z_{t+h}-A^z_{t+h}-(Y^z_{t}-A^z_{t}), \qq \pas
 \eeas
Since $Y^z_{t}-A^z_{t}-zB_{t} \in L^\infty(\cF_t)$, using  the
assumption (H5) one can check that
 \bea
 \label{Eeq0a}
 \cE\Big\{z(B_{t+h}-B_t)-\int_t^{t+h}g(s,z)ds\Big|\cF_t\Big\}=
 \cE\{Y^z_{t+h}-A^z_{t+h}|\cF_t\}-(Y^z_t-A^z_t)=0, \q\pas
 \eea
Therefore, applying (\ref{indep}) we have
 \beas
 hg(t,z)&=&\cE\Big\{z(B_{t+h}-B_t)-\int_t^{t+h}(g(s,z)-g(t,z))ds\Big|\cF_t\Big\}\\
 &=& \cE[z(B_{t+h}-B_t)|\cF_t]+v(t,h) =\cE[z(B_{t+h}-B_t)]+v(t,h),
 \eeas
where
 \beas
 v(t,h) &\dfnn& \cE\Big\{z(B_{t+h}-B_t)-\int_t^{t+h}(g(s,z)-g(t,z))ds\Big|\cF_t\Big\}
 -{\cal E}[z(B_{t+h}-B_t)|\cF_t]\\
    &=&\cE\Big\{zB_{t+h}-\int_t^{t+h}(g(s,z)-g(t,z))ds\Big|\cF_t\Big\}
    -\cE[zB_{t+h}|\cF_t].
 \eeas
Now, applying $L^p$-domination (\ref{pdom}) for the
$\cF$-expectation $\cE$ with $p=p\big(2\| \int_0^T|g(s,z)|ds
\|_\infty ,|z|\big)$, we obtain that
 \beas
E \Big\{\frac{1}{h^p}|v(t,h)|^p\Big\}&=& \|v(t,h) \|^p_p \leq 3^p
\Big\| \int_t^{t+h}|g(s,z)-g(t,z)|ds \Big\|^p_p\\
&\le & 3^p E\Big\{\frac{1}{h}\int_t^{t+h} |g(s,z)-g(t,z)|
ds\Big\}^p.
 \eeas
Since $z\in\hR^d$ is fixed, thus by the Lebesgue differentiation
theorem, $P$-almost surely one has
 \beas
 \frac{1}{h}\int_t^{t+h} |g(s,z)-g(t,z)| ds\rightarrow 0,\q \mbox{
 for a.e. }t\in [0,T].
 \eeas
The Dominated Convergence Theorem then implies that
 \beas
 E\Big\{\int_0^T
 \Big[\frac{1}{h}|v(t,h)|\Big]^p dt\Big\}\le 3^p E\Big\{\int_0^T
\Big[\frac{1}{h}\int_t^{t+h} |g(s,z)-g(t,z)| ds\Big]^p dt\Big\}
\to
 0.
 \eeas
In other words, we have proved that $v(t,h)=o(h)$ in
$\cH^p_\bF([0,T])$. Thus
 \beas g(t,z)=\lim_{h \to 0}
 \frac{\cE[z(B_{t+h}-B_t)]}{h},\q~~ \pas,
 \eeas
and it follows that $g$ is deterministic.

Now assume that $\cE$ also satisfies (H4), we show that $g$ is
local Lipschitz continuous.
To see this, taking $t+h=T$ in (\ref{Eeq0a}) and applying (H5) with
$\tilde{t}=\t=T$ we obtain that $\cE\big\{zB_T-\int_t^T
g(s,z)ds\big|\cF_t\big\}=zB_t$, $\pas$ Since $g$ is deterministic,
this implies that $ \int_t^T g(s,z)ds=\cE\{zB_T|\cF_t\}-zB_t$.
Similarly, one has $\int_t^T g(s,z')ds=\cE\{z'B_T|\cF_t\}-z'B_t$.
Combining, we have
 \beas \int_t^T [g(s,z')-g(s,z)]ds &=& \cE\{z
B_T|\cF_t\}-\cE\{z'B_T|\cF_t\}-(z-z')B_t\\
&\le &\cE^{z,z'}_\m\{(z-z')B_T|\cF_t\}-(z-z')B_t. \eeas
Note that for $g^{z,z'}_\m(v)\dfnn \m(1+|z|+|z'|)|v|$, one has
$$
\cE^{z,z'}_\m\{(z-z')B_T|\cF_t\}=(z-z')B_t+\int_t^T\m(1+|z|+|z'|)|z-z'|ds.
$$
We deduce that
$$\int_t^T [g(s,z')-g(s,z)]ds \le\int_t^T\m(1+|z|+|z'|)|z-z'|ds.
$$
Replacing $T$ by an arbitrary $t'\in(0,T]$ in the above, we can then
deduce that for any $t'\in(0,T]$, it holds that
$$
|g(t',z)-g(t',z')|\le \mu(1+|z|+|z'|)|z-z'|,
$$
proving the local Lipschitz property of $g$. \qed


The main result of this paper is the following representation
theorem.

\begin{thm} Assume that $\cE$ is a regular quadratic $\cF$-expectation that satisfies (H4)-(H6).
Then, there exists a local Lipschitz continuous function $g(t,z): [0,T]\times \hR^d \mapsto \hR$
such that for any $z \in \hR^d$,
 \bea
\label{gcontrol}
 g_1(t,\o,z) \leq g(t,z) \leq g_2(t,\o,z),\qq \dtp,
 \eea
and that for any $\xi \in L^\infty(\cF_T)$, it holds $P$-a.s. that
\beas \cE[\xi|\cF_t] = \cE^g[\xi|\cF_t],\qq \fa t \in [0,T]. \eeas
\end{thm}

{\it Proof.} Let $g$ be the random field defined in (\ref{repg}). We
know from Lemma \ref{dll} that $g$ is deterministic and local
Lipschitz continuous. Then (\ref{gcontrol}) follows from (\ref{e01})
and we see that $g\big|_{z=0}=0$. For any $\xi \in L^\infty(\cF_T)$,
we can apply the result of \cite[Theorem 2.3]{Ko} to conclude that
the BSDE($\xi,g$) admits a solution $(\hat{Y},\hat{Z}) \in
\hC^\infty_\bF([0,T])\times \cH^2_\bF([0,T];\hR^d)$. Furthermore, by
virtue of (\ref{locallip}), it follows from \cite{MM} (or
\cite{HIM}) that the solution is unique. (We remark that the result
of \cite{Ko} cannot be applied here since $g$ is not necessarily
differentiable).
Let $\{\P^n\}_{n \in
\hN}$ be a sequence of simple processes that approximates $\hat{Z}$
in $\cH^2_\bF([0,T];\hR^d)$. Then it holds that $\underset{t \in
[0,T]}{\sup}\Big|\int_0^t (\P^n_s-\hat{Z}_s) dB_s\Big|\to 0$ in
$L^2(\cF_T)$, thanks to the Burkholder-Davis-Gundy inequality.
Applying  \cite[Lemma 2.5]{Ko} we can find a subsequence of
$\{\P^n\}_{n \in \hN}$, still denoted by $\{\P^n\}_{n \in \hN}$,
such that
 \bea\label{Rconv} \P^n_t \to
 \hat{Z}_t,~~ \dtp  \q \mbox{and} \q \underset{t \in
 [0,T]}{\sup}\Big|\int_0^t (\P^n_s-\hat{Z}_s) dB_s\Big|\to 0,~~ \pas
 \eea
with $\underset{n \in \hN}{\sup}|\P^n_t| \in \cH^2_\bF([0,T])$ and
$\underset{n \in \hN}{\sup}\underset{t \in
[0,T]}{\sup}\Big|\int_0^t (\P^n_s-\hat{Z}_s) dB_s\Big| \in
L^2(\cF_T)$.
We define stopping times
 \bea\label{si_k}
 \si_k \dfnn \inf\Big\{t \in [0,T]: \int_0^t \underset{n \in \hN}{\sup}|\P^n_s|^2 ds
 +\underset{n \in \hN}{\sup}\underset{s \in [0,t]}{\sup}\Big|\int_0^s
 \P^n_r dB_r\Big|>k \Big\}\land T, \q~~ \fa k \in \hN.
\eea
It is easy to see that $\si_k \nearrow T$, $\pas$

For any $z \in \hR^d$, $0 \le t < \tilde{t}\le T$ and $\t \in
\cM_{0,T}$, it follows from (\ref{Eeq0a}) and (H6) that
 \bea\label{KEY0}
 \cE\Big\{\int_t^{\tilde{t}} \b1_{\{s \le \t \}}
 \big[-g(s, z)ds+ z dB_s\big]\Big|\cF_t\Big\}=0,\qq  \pas
 \eea
Let $\P$ be any member of $\{\P^n\}_{n \in \hN}$. Without loss of
generality we assume that $\P$ is in the form of \beas
\P_t(t,\o)=\sum^m_{i=0} \sum^{n_i}_{j=1}z^i_j
\b1_{[s_i,s_{i+1})\times E^i_j}(t,w),\qq \fa (t,\o) \in [0,T] \times
\O, \eeas where $0=s_0 < s_1<\cdot \cdot \cdot < s_m <s_{m+1}=T$,
$\;~\{E^i_j\}^{n_i}_{j=1}$ is an $\cF_{s_i}$-measurable partition of
$\O$ for $i =0,1\cdot  \cdot \cdot, m$, and each $z^i_j \in \hR^d$.

Now fix $k \in \hN$, for any $t \in [0,T]$, there exist $\a \in
\{0,1\cdot \cdot \cdot m\}$ such that $t \in [s_\a,s_{\a+1})$.
By refining the partition if necessary we may assume that $t=s_\a$.
Since the quadratic $\cF$-expectation $\cE$ is ``translation
invariant" and satisfies ``zero-one law'', using (\ref{KEY0}) one
can show that $P$-a.s.
 \bea \label{KEY}
 &&  \cE\Big\{\int_t^T \b1_{\{s \le \si_k\}} \big[-g(s, \P_s)ds+ \P_s dB_s\big]
 \Big|\cF_t\Big\} \nonumber\\
 &=&\cE\Big\{\sum^m_{i=\a} \sum^{n_i}_{j=1} \b1_{E^i_j} \int_{s_i}^{s_{i+1}}
 \b1_{\{s \le \si_k\}}\big[- g(s,z^i_j)ds +
 z^i_jdB_s\big]\Big|\cF_t\Big\}\nonumber\\
 & =& \cE\Big\{\neg \sum^{m-1}_{i=\a}\sum^{n_i}_{j=1}
\b1_{E^i_j} \int_{s_i}^{s_{i+1}}\b1_{\{s
\le \si_k\}} \big[- g(s,z^i_j)ds + z^i_jdB_s\big]\nonumber\\
&& +\sum^{n_m}_{j=1}\b1_{E^m_j}\cE\big[\int_{s_m}^T \b1_{\{s \le
\si_k\}} \big[- g(s,z^m_j)ds+ z^m_jdB_s\big] \big|\cF_{s_m}\big]
\Big|\cF_t\Big\}\nonumber\\
& =&  \cE\Big\{\sum^{m-1}_{i=\a}\sum^{n_i}_{j=1}\neg \b1_{E^i_j}
\neg \int_{s_i}^{s_{i+1}} \neg \b1_{\{s \le \si_k\}}\neg \big[\neg-
\neg g(s,z^i_j)ds + z^i_jdB_s\big]\Big|\cF_t\Big\}
\nonumber \\
&& \cdot  \cdot \cdot \cdot  \cdot \cdot\nonumber\\
& =& \cE\Big\{\sum^{n_\a}_{j=1} \b1_{E^\a_j} \int_t^{s_{\a+1}}
\b1_{\{s \le \si_k\}} \big[- g(s,z^\a_j)ds +
z^\a_jdB_s\big]\Big|\cF_t\Big\}\nonumber \\
& =& \sum^{n_\a}_{j=1}\b1_{E^\a_j}\cE\Big\{ \int_t^{s_{\a+1}}
\b1_{\{s \le \si_k\}} \big[- g(s,z^\a_j)ds +
z^\a_jdB_s\big]\Big|\cF_t\Big\}=0.
 \eea
For any $k \in \hN$, since $g$ is continuous and has quadratic
growth in $z$, using (\ref{Rconv}) and applying Dominated
Convergence Theorem we deduce that $\int_t^T \b1_{\{s \le \si_k\}}
[-g(s,\P^n_s)ds+  \P^n_s dB_s]$ converges to $\int_t^T \b1_{\{s \le
\si_k\}}[-g(s,\hat{Z}_s)ds+ \hat{Z}_s dB_s]$ almost surely. We also
see from the definition of $\si_k$ (\ref{si_k}) that
$$\Big|\int_t^T
\b1_{\{s \le \si_k\}} \big[-g(s,\P^n_s)ds+ \P^n_s dB_s ] \Big|\le
\ell T+2(1+\ell)k, \q \pas, \q \fa n \in \hN.
$$

Let $K=\ell T+2(1+\ell)k$ and $p\dfnn p(K,0)$, applying
$L^p$-domination of $\cE$ and using (\ref{KEY}) for each $\P^n$ one
can then deduce that $\cE\big\{\int_t^T \b1_{\{s \le \si_k\}}
\big[-g(s, \hat{Z}_s)ds+ \hat{Z}_s dB_s\big]\big|\cF_t\big\}=0$,
$\pas$
 \if{0}
 Moreover, the ``translation invariance" of $\cE$ and Proposition
\ref{cem} imply that $P$-a.s. for any $t \in [0,T]$,
 \bea
 \label{key}
 \cE\Big\{\int_0^T \b1_{\{s \le \si_k\}} \big[-g(s, \hat{Z}_s)ds+ \hat{Z}_s dB_s\big]
 \Big|\cF_t\Big\}=\int_0^t \b1_{\{s \le \si_k\}} \big[-g(s, \hat{Z}_s)ds+ \hat{Z}_s dB_s
 \big].
 \eea

Furthermore, for any $t \in [0,T]$, The optional sampling theorem
(Theorem \ref{opts}) implies that
 \beas Y^z_{t
 \land \si_n}-A^z_{t \land \si_n}=\cE[Y^z_{\si_n}-A^z_{\si_n}|\cF_t]
 ,\qq \pas
 \eeas
It then follows from (\ref{e02}) that
 \beas -\int_t^T \b1_{\{s \leq
 \si_n \}} g(s,z)ds +\int_t^T \b1_{\{s \leq \si_n \}}zdB_s
 =Y^z_{\si_n}-A^z_{\si_n}+A^z_{t \land \si_n}-Y^z_{t \land
 \si_n},\q~~ \pas
 \eeas
Since $\cE$ is translation invariant, one has, $P$-a.s.
 \bea
 \label{Eeq0}
 \cE\Big\{\neg-\neg \int_t^T \neg \neg \b1_{\{s \leq \si_n \}}
 g(s,z)ds \neg + \neg\int_t^T \neg \neg \b1_{\{s \leq \si_n
 \}}zdB_s\Big|\cF_t\Big\}\neg
 =\neg \cE[Y^z_{\si_n}\neg-\neg A^z_{\si_n}|\cF_t]\neg+\neg A^z_{t
 \land \si_n}\neg-\neg Y^z_{t \land \si_n}\neg=\neg 0.
 \eea
Moreover, Proposition \ref{cem} implies that $P$-a.s. for any $t \in
[0,T]$,
 \bea
 \label{key}
 \cE\Big\{\neg-\neg\int_0^T\neg\neg \b1_{\{s
 \leq \si_n \}} g(s,z)ds \neg+\neg\int_0^T\neg\neg \b1_{\{s \leq
 \si_n \}}zdB_s\Big|\cF_t\Big\}\neg=\neg-\int_0^t\neg \b1_{\{s \leq
 \si_n \}} g(s,z)ds \neg+\neg\int_0^t\neg \b1_{\{s \leq \si_n
 \}}zdB_s.
 \eea


 \no [{\bf Shige and Ying: Here is the problem. We need to show that
(\ref{key}) also holds when $z$ is replaced by the process $\dis Z
\in BMO$ (or probably just $\underset{p>0}{\bigcap}\cM^p(\hR^d)$.)
In your previous works this was done by approximation (in $L^2$),
combined with the $L^2$ domination. Here we can only use the
$L^p$-domination. The problem is that as you can see from Definition
\ref{Lpdom} that the exponent $p$ is not universal(!). We guess that
this will require some localization procedure, so you might have
better idea? One conjecture is that we might need to strengthen the
$L^p$-domination assumption for $\cE$ (see Appendix), but this
requires some justification as well. What do you think? (Note: the
following argument is based on assumption that this is true.)}] \ms
$\cdot \cdot \cdot \cdot \cdot \cdot \cdot $ Finally, for any $t \in
[0,T]$, \fi
The ``translation invariance" of $\cE$ then implies that
 \beas
 \cE[\hat{Y}_{
 \si_k}|\cF_t]\neg=\neg\cE\Big\{\hat{Y}_{t \land \si_k}
 +\int_t^T \b1_{\{s \le \si_k\}} \big[-g(s, \hat{Z}_s)ds+ \hat{Z}_s dB_s\big]
\Big|\cF_t\Big\}\neg =\hat{Y}_{t \land \si_k},
 \qq \pas
\eeas
%
Letting $p\dfnn p(\|\hat{Y} \|_\infty,0)$ and applying Theorem
(\ref{opts}) as well as $L^p$-domination for $\cE$ again, we obtain
that
 \beas \big\|\hat{Y}_{t \land \si_k} -\cE[\xi|\cF_t] \big\|_p =
 \big\|\cE[\hat{Y}_{\si_k}|\cF_t]-\cE[\xi|\cF_t] \big\|_p\leq 3 \|
 \hat{Y}_{\si_k}-\xi \|_p.
 \eeas
Since $\si_k \nearrow T$, $P$-a.s. and $\hat{Y}$ is continuous,
$\hat{Y}_{\si_k}$ converges $P$-a.s. to $\xi$ and $\hat{Y}_{t \land
\si_k}$ converges $P$-a.s. to $\hat{Y}_t$. These two convergence are
even in $L^p$ sense, thanks to the Lebesgue Dominated Convergence
Theorem. Thus $\cE[\xi|\cF_t] = \hat{Y}_t=\cE^g[\xi|\cF_t]$,
$P$-a.s. The conclusion then follows from Proposition \ref{cem} and
the continuity of $\hat{Y}$. \qed

\end{document}

\bibitem{BarKar}
Barrieu, P. and El Karoui, N. (2004) {\it Optimal derivatives
design under dynamic risk measures}, Article in Mathematics of
Finance , Contemporary Mathematics (A.M.S. Proceedings).

\bibitem{BCHMP}
Briand, P., Coquet, F., Hu, Y., M\'emain, J., and Peng, S. (2000),
{\it A converse comparison theorem for BSDEs and related
properties of $g$-expectations}, Electron. Comm. Probab. {\bf 5},
101-117.

\bibitem{BH-06}
Briand, P. and Hu, Y., preprint, \emph{BSDE with quadratic growth
and unbounded terminal value}

\bibitem{BH-07}
Briand, P. and Hu, Y., preprint, \emph{Quadratic BSDEs with Convex
Generators,}

\bibitem{CP-98}
Chen, Z. and Peng, S. (1998), {\it A nonlinear Doob-Meyer type
Decomposition and its application}, SUT Journal of Mathematics,
{\bf 34}(2), 197-208.

\bibitem{CHMP}
Coquet, F., Hu, Y., M\'emain, J., and Peng, S. (2002), {\it
Filtration-consistent nonlinear expectations and related
$g$-expectations}, Probab. Theory Relat. Fields, {\bf 123}, 1-27.

\bibitem{Ka}
Kazamaki, N. (1994) \emph{Continuous Exponential Martingales and BMO,}
Lecture Notes in Mathematics {\bf 1579}, Springer Berlin.

\bibitem{Ko} Kobylanski, M. (2000) \emph{Backward Stochastic Differential Equations
and Partial Differential Equations with Quadratic Growth,} Ann. Probab. {\bf 28},
no. 2, 558-602.

\bibitem{LSM-97}
Lepeltier, M. and San Martin, J. (1997), {\em Backward stochastic differential
equations with nonLipschitz coefficients}, {\sl Stat. and Prob. Letters},
v.32, {\bf  4}, 425--430.

\bibitem{LSM}
Lepeltier, J.-P. and San Martin, J. (1998) \emph{Existence for BSDE with superlinear-quadratic
cofficient,} Stochastics Stochastics Rep. {\bf 63}, no. 3-4, 227-240

\bibitem{MY}
Ma, J. and Yao, S.,preprint, (2006) \emph{Quadratic $g$-Expectations
and $g$-Martingales,}

\bibitem{mybk}
Ma, J. and Yong, J. (1999), {\sl Forward-Backward Stochastic
Differential Equations and Their Applications}, Lecture Notes in
Math., {\bf 1702}, Springer.

\bibitem{Peng-97}
Peng. S. (1997), {\it BSDE and related $g$-expectation}, In Pitman
Research Notes in Mathematics Series, No. 364, {\sl Backward
Stochastic Differential Equations, N. El Karoui and L. Mazliak
eds.}, 141-159.

\bibitem{Pln} Peng, S. \emph{Nonlinear Expectations, Nonlinear Evaluations and
Risk Measures,} Stochastic Methods in Finance,  \emph{Lectures given at the
C.I.M.E.-E.M.S. Summer School held in Bressanone/Brixen, Italy, July 6-12, 2003.}
 Lecture Notes in Mathematics {\bf 1856}, 165-253, Springer-Verlag.

\bibitem{Pr-90}
Protter, P. (1990), {\sl Integration and Stochastic Differential
Equations}, Springer.

\bibitem{Rev}
Revuz, D. and Yor, M. (1991), {\sl Brownian Motion and Continuous
Martingales}, Springer.

We first show that $g$ is deterministic and (\ref{limg}) holds. To
this end,
 \if{0} we note that if $g$
satisfies (H1) and (H2), one can deduce in the similar way as in
(\ref{qBSDE}) and (\ref{qBSDE1}) that \beas \big\{\xi+\int_0^T
\z_s dB_s: \xi \in L^\infty(\cF_T), ~ \z \in
L^\infty_\bF([0,T];\hR^d) \big\} \subset Dom(\cE^{g_1})\cap
Dom(\cE^{g_2})\subset Dom(\cE). \eeas
 \fi

 \bea
 \label{Eeq0a}
 \cE\big\{zB_{\tilde{t}\land \t}\dneg-\dneg
 \int_t^{\tilde{t}} \neg \b1_{\{s \le \t\}}g(s,z)ds\big|\cF_t\big\}
 \neg =\neg \cE\big\{Y^z_{\tilde{t}\land \t}\neg-\neg
 A^z_{\tilde{t}\land \t}\big|\cF_t\big\}\neg -\neg (Y^z_{t\land
 \t}\neg -\neg zB_{t \land \t}\neg -\neg A^z_{t \land \t})\neg=
 \neg zB_{t \land \t}.
 \eea
Taking $\tilde{t}=t+h$, $\t=T$, and noting that $zB_t-\cE[zB_t]$,
$t \in [0,T]$ is a $\cE$-martingale one shows that
 \beas
 hg(t,z)&\tneg=&\tneg \cE\big\{z B_{t+h}\neg-\neg\int_t^{t+h}\neg \big[g(s,z)-g(t,z)ds\big]\big|\cF_t\big\}-z B_t\\
 &\tneg= &\tneg \cE\big\{z B_{t+h}\neg-\neg\int_t^{t+h}\neg \big[g(s,z)-g(t,z)ds\big]\big|\cF_t\big\}-\cE\{zB_{t+h}|\cF_t\}+\cE\{zB_{t+h}\}-\cE\{zB_t\}.
 \eeas
Let $ v(t,h) \dfnn
\cE\big\{zB_{t+h}-\int_t^{t+h}(g(s,z)-g(t,z))ds\big|\cF_t\big\}
-\cE\{zB_{t+h}|\cF_t\}$, applying $L^p$-domination (\ref{pdom})
for the $\cF$-expectation $\cE$ with $p=p\big(2\|
\int_0^T|g(s,z)|ds \|_\infty ,|z|\big)$, we obtain that \beas
E\big[\frac{1}{h^p}|v(t,h)|^p\big]\dneg=\dneg \|v(t,h) \|^p_p \leq
\neg 3^p \Big\|\neg \int_t^{t+h}\neg |g(s,z)\neg-\neg g(t,z)|ds
\Big\|^p_p \neg \le \neg 3^p E\Big[\frac{1}{h}\neg\int_t^{t+h}\neg
|g(s,z)\neg-\neg g(t,z)|ds\Big]^p. \eeas Since $z\in\hR^d$ is
fixed, thus by the Lebesgue differentiation theorem, $P$-almost
surely one has
 \beas
 \frac{1}{h}\int_t^{t+h} |g(s,z)-g(t,z)| ds\rightarrow 0,\q \mbox{
 for a.e. }t\in [0,T].
 \eeas
The Dominated Convergence Theorem then implies that
 \beas
 E\Big\{\int_0^T
 \Big[\frac{1}{h}|v(t,h)|\Big]^p dt\Big\}\le E\Big\{\int_0^T
\Big[\frac{1}{h}\int_t^{t+h} |g(s,z)-g(t,z)| ds\Big]^p dt\Big\}
\to
 0.
 \eeas
In other words, we have proved that $v(t,h)=o(h)$ in
$\cH^p_\bF([0,T])$. Thus
 \beas g(t,z)=\lim_{h \to 0}
 \frac{\cE\{zB_{t+h}\}-\cE\{zB_t\}}{h},\q~~ \pas,
 \eeas
and it follows that $g$ is deterministic.
\section{Appendix}
\setcounter{equation}{0}

\begin{defn}
\label{strLpdom} A regular quadratic $\cF$-expectation $\cE$ is
said to satisfy the ``strong $L^p$-domination" if for any $K,R
>0$, there exist constants $p=p(K,R)>0$ and $C=C_R>0$ such that for any two
stopping times $0 \leq \t_1 \leq \t_2 \leq T$,
any $\xi_i \in L^\infty_{\t_i}$ with $\norm \xi_i \norm_\infty
\leq K,~i=1,2$ and any $\norm Z\norm_{BMO} \leq R$, it holds for
each $t \in [0,T]$ that
 \beas
 &&\big\|\cE\Big\{\xi_1+\int_0^{\t_1}Z_s
 dB_s\Big|\cF_t\Big\}-\int_0^{t \land \t_1} Z_s dB_s
 -\cE\Big\{\xi_2+\int_0^{\t_2}Z_s dB_s\Big|\cF_t\Big\}+\int_0^{t \land \t_2}
 Z_s dB_s \Big\|_p\\
 &&\qq\qq\qq \leq 3 \|\xi_1-\xi_2 \|_p+ C \|\t_1 - \t_2 \|_p.
 \eeas
\end{defn}

\section{Continuity of the Representation Generator}
\setcounter{equation}{0}

In the Lipschitz case, the proof of the continuity of the
representing generator depends crucially on the notion of
``domination", which is no longer valid in the quadratic case.
However, without such a control the proof of the continuity seems to
be impossible. Therefore, our first task is to find a reasonable
replacement of the concept of ``domination" of a $\cF$-consitent
nonlinear expectation.

We begin by considering a quadratic $g$-expectation. Suppose that
$g$ is deterministic function satisfying (H1) and (H2). It is
readily seen that $g$ must satisfy a ``local Lipschitz property":
\bea\label{locallip} |g(t,z)-g(t,z')|\le \m(1+|z|+|z'|)|z-z'|,
\qq\forall z,z'\in\hR^d. \eea Based on this observation, let us
denote, for each $z,z'\in\hR^d$, a function \bea\label{gmu}
g^{z,z'}_\m(v)\dfnn \m(1+|z|+|z'|)|v|, \qq \forall v\in\hR^d, \eea
and we denote the $g^{z,z'}_\m$-expectation to be
$\cE^{z,z'}_\m(\cd)$. It is worth noting that $\cE^{z,z'}_\m(\cd)$
is a Lipschitz $g$-expectation.

Now, again for fixed $z\in\hR^d$, consider the process
$Y^{g,z}_t\dfnn \cE^g\{zB_T|\cF_t\}$, $t\ge 0$. We should note that
this $g$-expectation is defined in a generalized sense, because the
terminal value $zB_T$ is not bounded, but rather one that has
exponential moment (e.g., $E\{e^{zB_T}\}=e^{\frac12
|z|^2T}<\infty$.)

\begin{lem}
\label{Ygzform} Suppose that $g$ is a deterministic function
satisfying (H1) and (H2). Then $Y^{g,z}_t$ has the following
explicit expression: \bea \label{Ygz} Y^{g,z}_t=zB_t+\int_t^T
g(s,z)ds, \qq t\in[0,T]. \eea
\end{lem}

{\it Proof.} Note that the main difficulty of this result comes from
the unboundedness of $B_T$, and the quadratic growth of $g$. To
overcome this difficulty we recall the definition of
$\cE^g\{zB_T|\cF_t\}$. That is, $Y^{g,z}_t$ is the monotone limit of
the sequence $\{Y^{g,z,n}_t\}_{n \ge 1}$, where each $Y^{g,z,n}$ is
the unique solution to the BSDE:
\bea\label{Ygzn} Y^{g,z,n}_t=zB_{T}+\int_{t}^{T}
g_n(s,Z^{g,z,n}_s)ds-\int_{t}^{T}Z^{g,z,n}_sdB_s, \qq t\in[0,T],
\eea
where $g_n$'s are Lipschitz in $z$ and $g_n\ua g$. Note that if $g$
is deterministic, then so are $g_n$'s. But for the deterministic
function $g_n$, the process
$$Y_t\dfnn zB_t+\int_t^T g_n(s,z)ds, \qq t\ge 0
$$
is an adapted solution to the BSDE (\ref{Ygzn}). Thus by uniqueness
$Y^{g,z,n}_t\equiv Y_t=zB_t+\int_t^T g_n(s,z)ds$, $t\ge 0$, for all
$n$. The Lemma then follows by letting $n\to\infty$. \qed

Now, let us define
\bea\label{hatE}
\hat{\cE}^{z,z'}\{(z-z')B_T|\cF_t\}(=\hat\cE^{z,z'}_t?)\dfnn
\cE^g\{zB_T|\cF_t\}-\cE^g\{z'B_T|\cF_t\},\qq t\ge 0. \eea

\begin{lem}
\label{hatEmu} Assume that $g$ is a deterministic function. Then the
process $\xi_t\dfnn\hat\cE^{z,z'}\{(z-z')B_T|\cF_t\}$, $t\ge 0$ is a
$\cE^{z,z'}_\m$-submartingale.
\end{lem}

{\it Proof.} First note that Lemma \ref{Ygzform} implies
\bea \label{hatEbsde} \hat\cE^{z,z'}_t&=&Y^{g,z}_t-Y^{g,z'}_t\\
&=& (z-z')B_T+\int_t^T(g(s,z)-g(s,z'))ds-\int_t^T(z-z')dB_s.
\nonumber \eea
For any $s\le t$, define
\bea \label{tildeY} \tilde{Y}_s&\dfnn& \cE^{z,z'}_\m\{\hat
\cE^{z,z'}_t|\cF_s\}=\Big[(z-z')B_t+\int_t^T(g(r,z)-g(r,z'))dr\Big]\nonumber\\
&&+\int_s^t\m(1+|z|+|z'|)|\tilde Z_r|dr-\int_s^t\tilde Z_rdB_r. \eea
But again, since $g$ is deterministic, the BSDE (\ref{tildeY}) has a
solution $(\hat Y,\hat Z)$, where
$$\hat Y_s\dfnn
(z-z')B_s+\int_t^T(g(r,z)-g(r,z'))dr+\int_s^t\m(1+|z|+|z'|)|z-z'|dr,
$$
and $\hat Z\equiv z-z'$.  Thus, denoting $\d g(r)\dfnn
g(r,z)-g(r,z')$, we have
\beas\tilde Y_t=\hat Y_t &=& (z-z')B_t+\int_t^T\d
g(r)dr+\int^t_s\m(1+|z|+|z'|)(z-z')dr+\int_s^t(z-z')dB_r\\
&=&(z-z')B_s+\int_s^T\d g(r)dr+\int^t_s\{\m(1+|z|+|z'|)(z-z')-\d
g(r)\}dr \\
&\ge & (z-z')B_s+\int_s^T \d g(r)dr. \eeas
But from (\ref{hatEbsde}) we see that the right hand side above is
exactly $
\hat\cE^{z,z'}_s=\xi_s$. This, combined with (\ref{tildeY}), shows
that $\xi=\hat\cE^{z,z'}$ is an $\cE^{z,z'}_\mu$-submartingale. \qed

Finally, let us consider a general nonlinear quadratic expectation.
Assume that \bi \item[(i)] there exists a constant $\m>0$, such that
for any fixed $z, z'$, it holds that \bea \label{Edomin} \cE(z
B_T|\cF_t\}-\cE(z'B_T|\cF\}\le \cE^{z,z'}_\m\{(z-z')B_T|\cF_t\};
\eea

\item[(ii)]
there exists a {\it deterministic} function $g$ such that \bea
\label{gidentity} \cE\Big\{-\int_t^T g(s,z)ds+z(B_T-B_t)|\cF_t\}=0.
\eea \ei We show that this is enough to guarantee that $g$ satisfies
the local Lipschitz condition (\ref{locallip})!

\ss Indeed, since $g$ is deterministic, (\ref{gidentity}) implies
that
$$ \int_t^T g(s,z)ds=\cE(z(B_T-B_t)|\cF_t\}=\cE(zB_T|\cF_t\}-zB_t.
$$
Similarly, one has
$$\int_t^T g(s,z')ds=\cE(z'(B_T-B_t)|\cF_t\}=\cE(z'B_T|\cF_t\}-z'B_t.
$$
Combining, we have
 \beas \int_t^T [g(s,z')-g(s,z)]ds &=& \cE(z
B_T|\cF_t\}-\cE(z'B_T|\cF\}-(z-z')B_t\\
&\le &\cE^{z,z'}_\m\{(z-z')B_T|\cF\}-(z-z')B_t. \eeas
Note that for $g^{z,z'}_\m(v)\dfnn \m(1+|z|+\z'|)|v|$, one has
$$
\cE^{z,z'}_\m\{(z-z')B_T|\cF_t\}=(z-z')B_t+\int_t^T\m(1+|z|+|z'|)|z-z'|ds.
$$
We deduce that
$$\int_t^T [g(s,z')-g(s,z)]ds \le\int_t^T\m(1+|z|+|z'|)|z-z'|ds.
$$
Replacing $T$ by an arbitrary $t'\in(0,T]$ in the above, we can then
deduce that for any $t'\in(0,T]$, it holds that
$$
|g(t',z)-g(t',z')|\le \mu(1+|z|+|z'|)|z-z'|,
$$
proving the local Lipschitz property of $g$.

\begin{lem}
If we additionally assume that for any $z \in \hR^d$ and any $t \in
[0,T]$ \beas \cE[z(B_T-B_t)|{\cal F}_t]=\cE[z(B_T-B_t)],\q~~ \pas
\eeas then $g$ is deterministic.
\end{lem}
\begin{proof}
Since \restart \beas
\cE[z(B_{t+h}-B_t)-\int_t^{t+h}g(s,z)ds|\cF_t]=0, \eeas we have
\beas
hg(t,z)&=&\cE[z(B_{t+h}-B_t)-\int_t^{t+h}(g(s,z)-g(t,z))ds|\cF_t]\\
&=& \cE[z(B_{t+h}-B_t)|\cF_t]+v(t,h)\\
&=& \cE[z(B_{t+h}-B_t)]+v(t,h), \eeas where \beas
v(t,h) &\dfnn& \cE[z(B_{t+h}-B_t)-\int_t^{t+h}(g(s,z)-g(t,z))ds|\cF_t]-{\cal E}[z(B_{t+h}-B_t)|\cF_t]\\
    &=&\cE[zB_{t+h}-\int_t^{t+h}(g(s,z)-g(t,z))ds|\cF_t]-\cE[zB_{t+h}|\cF_t].
\eeas Let for $p=p\big(2\norm \int_0^T|g(s,z)|ds \norm_\infty
,|z|\big)$, applying $L^p$-domination for $\cE$ we obtain that \beas
 \norm v(t,h) \norm_p \leq 3 \norm \int_t^{t+h}|g(s,z)-g(t,z)|ds
 \norm_p.
\eeas Hence, one has \beas E \big[\frac{1}{h^p}|v(t,h)|^p\big] \le
3^p E \Big[\frac{1}{h^p}\big(\int_t^{t+h} |g(s,z)-g(t,z)|
ds\big)^p\Big]. \eeas Lebesgue's theorem then implies that $P$-a.s.
\beas \frac{1}{h}\int_t^{t+h} |g(s,z)-g(t,z)| ds\rightarrow 0,\q
\mbox{ for a.e. }t\in [0,T], \eeas or $P$-a.s.  \beas \int_0^T
(\frac{1}{h}\int_t^{t+h} |g(s,z)-g(t,z)| ds)^p dt\rightarrow 0.
\eeas Thus, we can deduce that \beas E\big[\int_0^T
\big(\frac{1}{h}|v(t,h)|\big)^p dt\big]\le E\Big[\int_0^T
\big(\frac{1}{h}\int_t^{t+h} |g(s,z)-g(t,z)| ds\big)^p dt\Big] \to
0. \eeas Applying Lebesgue Convergence theorem, we see that
$v(t,h)=o(h)$ in $\cH^p_\cF$. Finally, \beas g(t,z)=\lim_{h \to 0}
\frac{\cE[z(B_{t+h}-B_t)]}{h},\q~~ \pas, \eeas which means that $g$
is deterministic.
\end{proof}